\documentclass[10pt]{article}
\usepackage{amsmath, amsfonts, amsthm}
\usepackage{amssymb,mathrsfs, verbatim}
\usepackage{graphicx,amstext}
\usepackage{amsmath}
\usepackage{mathtools}
\usepackage{amssymb}
\usepackage{latexsym}
\usepackage{ mathrsfs }
\usepackage{enumerate}
\usepackage[USenglish]{babel}

\usepackage[latin1]{inputenc}
\newcommand{\clg}[1]{{\mathcal{#1}}}

\newcommand{\ol}[1]{{\overline{#1}}}

\newcommand{\R}{\mathbb R}

\newcommand{\N}{\mathbb N}

\newcommand{\vp}{\varphi}

\newcommand \loc    {\text{loc}}
\newcommand{\ess}{\,{\rm ess}}

\newcommand{\dive}{{\rm div}}

\newcommand{\dist}{\,{\rm dist}}

\numberwithin{equation}{section}

\newtheorem{theorem}{Theorem}[section]
\newtheorem{proposition}{Proposition}[section]
\newtheorem{remark}{Remark}[section]
\newtheorem{lemma}{Lemma}[section]
\newtheorem{corollary}{Corollary}[section]
\newtheorem{definition}{Definition}[section]

\newcommand{\intav}[1]{\mathchoice {\mathop{\vrule width 6pt height 3 pt depth  -2.5pt
\kern -8pt \intop}\nolimits_{\kern -6pt#1}} {\mathop{\vrule width
5pt height 3  pt depth -2.6pt \kern -6pt \intop}\nolimits_{#1}}
{\mathop{\vrule width 5pt height 3 pt depth -2.6pt \kern -6pt
\intop}\nolimits_{#1}} {\mathop{\vrule width 5pt height 3 pt depth
-2.6pt \kern -6pt \intop}\nolimits_{#1}}}
 
\begin{document}
\title{Initial mixed-boundary value problem 
\\
for anisotropic fractional degenerate 
\\
parabolic equations}
\author{Gerardo Huaroto $^1$, Wladimir Neves $^2$}
\date{}

\maketitle

\footnotetext[1]{Instituto de Matem\'atica, Universidade Federal de Alagoas, Maceio, Alagoas, Brazil. E-mail: {\sl gerardo.cardenas@im.ufal.br}}

\footnotetext[2]{Instituto de Matem\'atica, Universidade Federal
do Rio de Janeiro, Rio de Janeiro, Rio de Janeiro, Brazil. E-mail: {\sl wladimir@im.ufrj.br.}
}
\begin{abstract}
We consider an initial mixed-boundary value problem for anisotropic fractional type degenerate parabolic equations
posed in bounded domains.
Namely, we consider that the boundary of the domain splits into two parts.
In one of them, we impose a Dirichlet boundary condition and in the another one a
Neumann condition. Under this mixed-boundary condition, 
we show the existence of solutions for
measurable and bounded non-negative initial data.
The nonlocal anisotropic diffusion effect relies on an inverse
of a $s-$fractional type elliptic operator, and 
the solvability is proved for any $s \in (0,1)$.
\end{abstract}

\maketitle

\section{Introduction}
\label{INTRO}
We are concerned in this paper with an initial mixed-boundary value problem for a class of anisotropic fractional 
type degenerate parabolic equations. To this end, 
let $\Omega\subset \R^n$ be a bounded open set with  smooth ($C^2$) boundary $\Gamma$, and 
denote by $\nu$ the outward unit normal vector field on it.
We assume that $\Gamma$ is divided into two parts $\Gamma_0$, $\Gamma_1$. 
Then, we consider the following initial mixed-boundary value problem
\begin{equation}
\label{FTPME} 
     \left \{
       \begin{aligned}
       \partial_t u + \dive  \mathbf{q}&= 0 \quad \text{ in $\Omega_T$},
       \\
       u|_{\{t=0\}}&= u_0 \quad \text{in $\Omega$},
       \\
       u &= 0 \quad \hspace{5pt} \text{on $(0,T) \times \Gamma_0$},
       \\
      \mathbf{q} \cdot \nu &= 0 \quad \hspace{5pt} \text{on $(0,T) \times \Gamma_1$},
\end{aligned}
\right.
\end{equation}
where $\Omega_T= (0,T) \times \Omega$, for any real number $T> 0$,
$u(t,x)$ is a real function, which could be interpreted as a density (concentration, population, etc.) or the thermodynamic temperature, 
$\mathbf{q}= - u \ A(x) \, \nabla \mathcal{K}_s u$ is the diffusive fractional flux, and $\mathcal{K}_s$ is the 
inverse of the $s$-fractional elliptic operator $\mathcal{L}^{s}_{\mathcal{B}}$, $(0< \!s < \!1)$, see Section \ref{DNSFO}. 
The matrix $A(x)= (a_{ij}(x))_{n\times n}$
is assumed symmetric and satisfies 
\begin{equation}
\label{REGA}
     a_{ij} \in C(\overline{\Omega}) \cap C^{0,1}_{\loc}(\Omega), \quad (i,j= 1,\ldots,n), 
\end{equation}
\begin{equation}
\label{unifellip}
  \sum_{i,j=1}^na_{ij}(x)\xi_i\xi_j \geq \Lambda_1 |\xi|^2, 
\end{equation}
for all $\xi\in\R^n$ and each $x \in \Omega$, for some ellipticity constant 
$\Lambda_1> 0$.
Moreover, the initial data $u_0 \in L^\infty(\Omega)$ is a non-negative given function, 
and we consider homogeneous 
Dirichlet and Neumann 
boundary conditions, respectively on $\Gamma_0$, $\Gamma_1$.   
This assumption, the mixed-boundary condition,
brings some difficulties which are discussed
through this paper, see for instance Section \ref{BOCOMEN}.  

\medskip
The diffusive non-local flux $\mathbf{q}$ in the initial mixed-boundary value problem 
\eqref{FTPME} is motivated by the so-called 
General Fractional Fick's law 
$$
  \mathbf{q}(x,u):= -\kappa(x,u) \ \nabla \clg{F} u
$$
provided $\kappa(\cdot,u)$ is positive (non-negative in general) defined, 
where $\clg{F}$  is the 
inverse of a fractional elliptic operator.
The first attempted is to consider 
$$
  \mathbf{q}(x,u):= - g(u) A(x) \nabla \mathcal{K}_s u
$$
with $g(u)= u$ or $g(u)= u (1 - u)$, which from the maximum principle
ensures that, $\kappa(\cdot,u)$ is non-negative defined. 
For the second case, $g(u)= u (1 - u)$, it should be also assumed that,
$0 \leq u_0 \leq 1$, but we leave this option to future work (see \cite{GHWN3}).   
Moreover, the assumption here $\kappa(x,u)= u A(x)$ turns clear that, 
the coefficients $(a_{ij})$, $(i,j= 1,\ldots,n)$ describe the anisotropic and
the heterogeneous nature of the medium. This is very important 
to a great many physical theories, for instance, let us mention applications in 
physical-chemical reactions and biological processes. 
Although, it is essential to mention that, 
in another context of porous media diffusion model, 
Caffarelli, Vazquez \cite{Caffa} introduced for the first time the model \eqref{FTPME} for 
a given fractional potential pressure law, that is to say, they considered 
$\mathbf{q}(u)= - u \nabla \clg{K} u$, where $\clg{K}$ is the inverse of the $s-$fractional Laplacian in $\R^n$. 
Hence that paper established a Fractional Darcy's law and under some conditions, 
they proved existence of weak (non-negative) solutions for the Cauchy problem.

\medskip
Concerning the elliptic linear operator $\mathcal{L}u:= -\dive(A(x)\nabla u)$,
which is the block building for the construction of the fractional operator 
$\mathcal{L}\,^{s}_{\mathcal{B}}$, we were motivated by the paper of Caffarelli, 
Stinga \cite{LCPRS}. In that paper the authors reproduce
Caccioppoli type estimates (for the Dirichlet and also Neumann boundary conditions),
which alloy them to develop the interior and boundary regularity theory
depending on the smoothness of the matrix $A(x)$ and the source terms. 
Albeit, we should mention that, different from that paper, 
here we are focused in the minimal regularity for the matrix $A(x)$, such 
that, the eigenfunctions $\{\varphi_k\}$ of the problem \eqref{EVP} have the enough
regularity to define conveniently the operator $\mathcal{K}_s$, and also 
to give a sense of the Neumann boundary condition on $\Gamma_1$, 
that is to say, for each function $\gamma \in H^1_0(0,T; H^1_{\Gamma_0}(\Omega))$ 
\begin{equation}
\ess\lim_{\tau \to 0^+}\int_0^T \!\! \int_{\Gamma_1} \mathbf{q}\left(\Psi_\tau(r), u(t,\Psi_\tau(r)) \right)\cdot \nu_\tau(\Psi_\tau(r)) \ \gamma(t,r) \,dr\,dt= 0,\nonumber
\end{equation}
where $\Psi_\tau(r):=r-\tau\nu(r)$, and $\nu_\tau$ is the unit outward normal field on $\Psi_\tau(\Gamma)$,
see Appendix. Recall that $A(x)$ is (uniformly) 
continuous up to the boundary, therefore it is bounded in $\Omega$ and its restriction 
on $\Gamma$ makes sense. This is also important to the $\clg{L}_{\clg{B}}$ operator's domain definition,
see equation \eqref{DOLB}. 
Moreover, due to the regularity of the matriz $A $ in $C^{0,1}_{\loc}(\Omega)$,
the eigenfunctions $\vp_k \in H^2(\Omega')$, for all $k \geq 1$ and every 
$\Omega'$ compactly contained in $\Omega$, see L. Ambrosio et al. \cite{ACM}. 
We remark that, it is not possible to ensure $H^2(\Omega)$ regularity even if the diffusive matriz $A$ has $C^{0,1}(\Omega)$ 
smoothness. Indeed, we are considering mixed-boundary conditions and hence Nirenberg's type methods do not apply, since
$\vp_k= 0$ on $\Gamma_0$ but not necessarily zero on $\Gamma_1$. 

\medskip
Since the paper \cite{Caffa}, there exists a considerable list of important correlated results, to mention a few
\cite{BIK}, \cite{CSV}, \cite{CV2}, \cite{I}, \cite{LMS}, \cite{STV2}, \cite{STV3}, \cite{STV5}. 
In particular, along the same problem 
Caffarelli, Soria and Vazquez establish the H\"older regularity of such weak 
solutions for the case $s \neq1/2$ in \cite{CSV}, and the case  $s= 1/2$ has been proved by Caffarelli, Vazquez
in \cite{CV2}. All of these above cited papers are posed in 
$\R^n$. On the other hand, the authors considered in \cite{Hua} again $\mathbf{q}(u)= - u \nabla \clg{K} u$, 
but now in the context of heat equation (Fractional Fourier law), and it was 
considered homogeneous Dirichlet boundary condition.
Thus the problem were posed in a bounded open subset of $\R^n$.   
One of the main task of that paper was how the boundary condition should be assumed, 
and it was important to deal with traces at the boundary for any $s \in (0,1)$. The problem 
here has new different difficulties, and a different context. In this way we consider another formulation different from 
that one presented in \cite{Hua}. Indeed,
an important pragmatism concerning the mixed-boundary conditions is that, 
the (homogeneous) Dirichlet boundary conditions are taken into account in the test functions, 
and the Neumann boundary conditions are taken into account in the linear form due to boundary 
integrals. Hence we follow this strategy and address the reader to 
Section \ref{SOLVAB}, where the main ideas are well-explained and also Section \ref{MainResult},
where it is shown the solvability of the initial mixed-boundary value problem \eqref{FTPME}. 

\medskip
Finally, we would like to stress that the uniqueness property is not established in this paper.
First, let us remark that, no uniqueness result has been proven
even for the $\R^n$ case with $\mathbf{q}(u)= - u \nabla \clg{K} u$. 
Moreover, along the same model we address the reader to Serfaty, V\'azquez \cite{SSJV} (and references therein),  
where is constructed a counterexample to comparison of densities,
see Section 6.5 (Lack of comparison principle). Hence we may consider 
a selection principle (or admissibility criteria) in order to attack the 
issue of uniqueness for \eqref{FTPME}.

\subsection{Functional Space}

From now on by $\Omega$ we denote a bounded open set in $\R^n$ with smooth ($C^2$) boundary $\Gamma$. 
We assume that $\Gamma= \Gamma_0 \cup \Gamma_1$, $\Gamma_0$ is a closed set and $\mathcal{H}^{n-1}(\Gamma_0)> 0$, where
$\mathcal{H}^\theta$ is the usual $\theta-$Hausdorff measure.
Moreover, $\Gamma_0\cap\overline{\Gamma}_1$ is a submanifold of codimension greater than 1.
%
Then, we define
$$
   H^1_{\Gamma_0}(\Omega):= \left\lbrace v\in H^1(\Omega): \text{$v= 0$ on $\Gamma_0$ in the sense of trace} \right\rbrace,
$$
endowed with the norm
\begin{equation}
\label{normH1Gama0}
    \| v \|_{H^1_{\Gamma_0}(\Omega)}:=\left(\int_\Omega|\nabla v(x)|^2 \ dx\right)^{1/2},
    \quad \text{for each $v\in H^1_{\Gamma_0}(\Omega)$}.
\end{equation}
Since the trace is a continuous operator, we have that $H^1_{\Gamma_0}(\Omega)$ is a Hilbert space 
with the norm $\|\cdot\|_{H^1(\Omega)}$, which is equivalent to \eqref{normH1Gama0}.
Moreover, we define the set 
\begin{equation}
\label{DEFTESTFUNC}
C_{\Gamma_0}^{\infty}(\overline{\Omega}):=\left\lbrace  v \in C^{\infty}(\overline{\Omega});\, v= 0 \mbox{ on } \Gamma_0\right\rbrace,
\end{equation}
which is dense in $H^1_{\Gamma_0}(\Omega)$.


\bigskip
Now, we follow Lions, Magenes \cite{LionsMagenes} for the definition
of the spaces $H^s(\Omega)$, with $s \in (0,1)$. Indeed,  
by interpolation between $H^1(\Omega)$ and $L^2(\Omega)$, we have
$$
   H^s(\Omega)=[H^1(\Omega),L^2(\Omega)]_{1-s}.
$$
According to this definition, this space is a Hilbert space with the natural norm given by the interpolation. Moreover, we can define the space $H^s_0(\Omega)$ by
$$
H^s_0(\Omega)=\overline{C^{\infty}_c(\Omega)}^{\|\cdot\|_{H^s(\Omega)}}.
$$

Since $\Omega$ has regular boundary, the set $H^s_0(\Omega)$
could be written as an interpolation (see Theorem 11.6 of \cite{LionsMagenes}), 
$$
H^s_0(\Omega)=[H^1_0(\Omega),L^2(\Omega)]_{1-s},
$$
for each $s\in (0,1) \setminus \{1/2\}$. The particular case $s = 1/2$ generates the so called 
Lions-Magenes space $H^{1/2}_{00}(\Omega)$, which is defined by
$$
    H^{1/2}_{00}(\Omega):= [H^1_0(\Omega),L^2(\Omega)]_{1/2},
$$
which has the following characterization
$$
H^{1/2}_{00}(\Omega)=\left\lbrace u\in H^{1/2}(\Omega);\int_{\Omega}\dfrac{u(x)^2}{\dist(x,\Gamma)}dx<\infty\right\rbrace.
$$

Furthermore, we define the space $H^s_{\Gamma_0}(\Omega)$ by
$$
H^s_{\Gamma_0}(\Omega)=\mbox{ closure of $C^{\infty}_{\Gamma_0}(\bar{\Omega})$ in $H^s(\Omega)$}.
$$
In particular, for $0< s \leq 1/2$ and since $\Gamma$ is Lipschitz, we have 
$H^s_{\Gamma_0}(\Omega)=H^s(\Omega)$, which is due to the fact that $C^{\infty}_0(\Omega)$ is dense in $H^s(\Omega)$ 
(see \cite{LionsMagenes} Theorem 11.1).
On the other hand, if $1/2<s<1$ and $\Gamma$ is Lipschitz, then the spaces $H^s_{\Gamma_0}(\Omega)$ have a characterization via Trace operator (Theorem 9.4 \cite{LionsMagenes}), hence
\begin{equation}
   H_{\Gamma_0}^s(\Omega) \equiv \{ u \in H^s(\Omega): \text{$u= 0$ on $\Gamma_0$ in the sense of trace}\}. \label{trace}
\end{equation}
The proof is based in similar arguments considered in Theorem 11.5
 \cite{LionsMagenes}.

\medskip
Finally, since $\Omega$ has Lipschitz boundary,
there exists an equivalent definition given via interpolation. Indeed, due to 
$H^1_0(\Omega)\subset H^1_{\Gamma_0}(\Omega)\subset H^1(\Omega)$, it follows that, 
for all $s \in (0,1)$
$$
[H^1_{0}(\Omega),L^2(\Omega)]_{1-s}\subset[H^1_{\Gamma_0}(\Omega),L^2(\Omega)]_{1-s}\subset[H^1(\Omega),L^2(\Omega)]_{1-s}.
$$
Therefore, we have
\begin{equation}
\begin{aligned}
H^s_0(\Omega)&\subset[H^1_{\Gamma_0}(\Omega),L^2(\Omega)]_{1-s}\subset H^s(\Omega),\quad s\in(0,1) \setminus \{1/2\}
\\[5pt]
H^{1/2}_{00}(\Omega)&\subset[H^1_{\Gamma_0}(\Omega),L^2(\Omega)]_{1/2}\subset H^{1/2}(\Omega),\quad s=1/2.
\end{aligned}\label{Eq_H_inclu}
\end{equation}
In particular, when $0<s<1/2$ we obtain
$$
[H^1_{\Gamma_0}(\Omega),L^2(\Omega)]_{1-s}=H^s(\Omega).
$$
On the other hand, using the idea of Theorem 11.6 \cite{LionsMagenes} we may obtain
$$
[H^1_{\Gamma_0}(\Omega),L^2(\Omega)]_{1-s}=H^s_{\Gamma_0}(\Omega),\quad \text{for all $s\in(1/2,1)$}.
$$

\section{Dirichlet-Neumann Spectral Fractional \\ Elliptic Operators}
\label{DNSFO}

In this section, we study some results of Dirichlet-Neumann spectral fractional elliptic operators.  
We mainly provide the proofs of the new results, in particular we stress Proposition \ref{proKuHu}. 
One can refer to \cite{MSireVazquez}, \cite{LCPRS}, and \cite{Hua} for an introduction.

\medskip
We are mostly interested in fractional powers of a strictly positive self-adjoint 
operators defined in a domain, which is dense in a (separable) Hilbert space. Therefore,
we are going to consider the linear operator $\mathcal{L}u= -\dive(A(x)\nabla u)$ equipped with 
homogeneous mixed Dirichlet-Neumann boundary data, that is to say $\mathcal{B}(u)=0$ on $\Gamma$, 
where the boundary operator $\mathcal{B}$ is defined as follows
\begin{equation}
\label{DNdata} 
\mathcal{B}(u)=\left \{
       \begin{aligned}
       u  \quad &\text{on $\Gamma_0$},
       \\
       (A \, \nabla u) \cdot \nu \quad &\text{on $\Gamma_1$},
\end{aligned}
\right.
\end{equation}
where $A(x)$ is the symmetric matrix satisfying 
\eqref{REGA} and \eqref{unifellip}.

\medskip
For conveniency, let us denote by $\mathcal{L}_{\mathcal{B}}$ the operator $\mathcal{L}$ 
subject to Dirichlet-Neumann boundary condition given by \eqref{DNdata}. 
Observe that $\mathcal{L}_{\mathcal{B}}$ is nonnegative and selfadjoint in
\begin{equation}
\label{DOLB}
D(\mathcal{L}_{\mathcal{B}}):= \left\lbrace u\in H^1(\Omega):\dive(A\nabla u)\in L^2(\Omega),\mbox{with $\mathcal{B}(u)= 0$ on $\Gamma$} \right\rbrace.
\end{equation}
 Therefore, by the spectral theory, there exists a complete orthonormal 
basis $\{\vp_k\}^{\infty}_{k=1}$ of $L^2(\Omega)$,
where 
$\vp_k$ satisfies
\begin{equation}
\label{EVP}
\left\lbrace
\begin{aligned}
\mathcal{L}\, \vp_k&=\lambda_k\vp_k,\quad\mbox{ in }\Omega,
\\[5pt]
\mathcal{B}(\vp_k)&= 0,\quad\quad\quad\mbox{ on }\Gamma.
\end{aligned}
\right.
\end{equation} 
It is easy to check that $\{\vp_k\}^{\infty}_{k=1}$ is also an orthogonal basis of $H^1_{\Gamma_0}(\Omega)$. 
Moreover, due to the regularity of the matriz $A(x)$ the eigenfunctions $\vp_k \in H^2(\Omega')$, for all $k \geq 1$ and every 
$\Omega'$ compactly contained in $\Omega$, see L. Ambrosio et al. \cite{ACM}. 


\medskip
For each $k\geq 1$, it follows that $\vp_k$ is an eigenfunction corresponding to $\lambda_k$, 
where one repeats each eigenvalue $\lambda_k$ according to its (finite) multiplicity
$$
    0< \lambda_1\leq\lambda_2\leq\lambda_3\leq \cdots \leq \lambda_k  \leq \cdots, \quad \text{$\lambda_k \rightarrow \infty$
     as $k \longrightarrow \infty$}.
$$
Then, we have 
$$
\begin{aligned}
   D(\,\mathcal{L}_{\mathcal{B}}\,)&= \{u \in L^2(\Omega); \ \sum_{k= 1}^\infty \lambda_k^2 \, |\langle u, \vp_k \rangle|^2 < \infty \}, 
   \\[5pt]
   \mathcal{L}_{\mathcal{B}} \ u&= \sum_{k= 1}^\infty \lambda_k \, \langle u, \vp_k \rangle \ \vp_k, 
   \quad \text{for each  $u \in D(\mathcal{L}_{\mathcal{B}})$}.
\end{aligned}
$$
Now, applying the functional calculus, we define for each $s >0$,
the following fractional elliptic operator
 $\mathcal{L}\,^s_{\mathcal{B}}$, given by
 $$
    \mathcal{L}\,^s_{\mathcal{B}}\, u:= \sum_{k= 1}^\infty \lambda^s_k \, \langle u, \vp_k \rangle \ \vp_k, 
$$ 
and it is well defined in the space of functions
\begin{equation}
\label{DEFFRALAP}
\begin{aligned}
    D\big(\,\mathcal{L}\,^s_{\mathcal{B}}\,\big)&= \Big\{ u \in L^2(\Omega) : \; \sum^{\infty}_{k=1}
    \lambda^{2 s}_k \, \vert \langle u, \vp_k \rangle \vert^2<+\infty  \Big\},
\end{aligned}
\end{equation}
which is a Hilbert space with the inner product
\begin{equation}
      \langle u, v \rangle_{s}:=\langle u, v\rangle+\int_{\Omega} \mathcal{L}\,^s_{\mathcal{B}}\, u(x) \; \mathcal{L}\,^s_{\mathcal{B}}\, v(x) \, dx. \nonumber
\end{equation}
In particular, the norm $\vert \cdot \vert_s$ is defined by
\begin{equation}
\label{eq:norm}
       \vert u \vert_s^2=\Vert u \Vert^2_{L^2(\Omega)}+\Vert \,\mathcal{L}\,^s_{\mathcal{B}}\, u \Vert^2_{L^2(\Omega)}.
\end{equation}

\medskip
Analogously, we can also define $\mathcal{L}\,^{-s}_{\mathcal{B}}: D\big(\mathcal{L}\,^{-s}_{\mathcal{B}}\big)\subset L^2(\Omega) \to L^2(\Omega)$. 
The next proposition give us the main properties of the operators defined above.
In particular, we observe that $D\big(\mathcal{L}\,^{-s}_{\mathcal{B}}\big)= L^2(\Omega)$.
\begin{proposition}
\label{THMCARAC}
Let $\Omega \subset \R^n$ be a bounded open set with Lipschitz boundary, $s \in (0,1)$, and consider 
the operators $\mathcal{L}\,^{s}_{\mathcal{B}}$, and $\mathcal{L}\,^{-s}_{\mathcal{B}}$. Then, we have:

\begin{enumerate}
\item[$(1)$] The operator   $\mathcal{L}\,^{s}_{\mathcal{B}}$ and $\mathcal{L}\,^{-s}_{\mathcal{B}}$ are self-adjoint.
Also $(\mathcal{L}\,^{s}_{\mathcal{B}})^{-1}= \mathcal{L}\,^{-s}_{\mathcal{B}}$.

\item[$(2)$] If $0 \leq s_{1}<s_{2} \leq 1,$ then
$$
D\left(\mathcal{L}\,^{s_2}_{\mathcal{B}}\right) \hookrightarrow D
\left(\mathcal{L}\,^{s_1}_{\mathcal{B}}\right),
\mbox{ and } D\left(\mathcal{L}\,^{s_2}_{\mathcal{B}}\right) \mbox{ is dense in } D\left(\mathcal{L}\,^{s_1}_{\mathcal{B}}\right).
$$

\item[$(3)$] For each $s,\,\sigma>0$ and $u\in D(\mathcal{L}\,^{s}_{\mathcal{B}})$ we have
$\mathcal{L}\,^{-\sigma}_{\mathcal{B}}\,u\in D(\mathcal{L}\,^{s+\sigma}_{\mathcal{B}})$.
\end{enumerate}
\end{proposition} 
\begin{proof}
The proof proceed analogously to the Proposition 2.1 in \cite{Hua} and hence we omit it.
\end{proof}

\medskip
Now, we state a Poincare's type inequality for the $\mathcal{L}\,^{s}_{\mathcal{B}}$, and an equivalent 
norm for $D\big(\,\mathcal{L}\,^{s}_{\mathcal{B}}\,\big)$.

\begin{corollary}[ Poincare's type inequality ]
\label{poincare}
Let $\Omega \subset \R^n$ be a bounded open set with Lipschitz boundary. Then for each $s> 0$, we have
$$
    \Vert u \Vert_{L^2(\Omega)} \leq \lambda_1^{-s} \ \Vert \mathcal{L}\,^{s}_{\mathcal{B}}\, u \Vert_{L^2(\Omega)}, 
    \quad \text{for all $u \in D\big(\,\mathcal{L}\,^{s}_{\mathcal{B}}\,\big)$}.
$$
Moreover, the norm defined in \eqref{eq:norm} and
\begin{equation}
\label{equi}
  \Vert u \Vert_s^2:= \int_{\Omega}\vert \mathcal{L}\,^{s}_{\mathcal{B}}\, u(x) \vert^2 \ dx
\end{equation}
are equivalent. 
\end{corollary}
\begin{remark}
As a consequence of the above results, we could consider the  inner product in 
$D\big(\mathcal{L}\,^{s}_{\mathcal{B}}\big)$, as follow
\begin{equation}
     \langle u, v \rangle_{s}= \int_{\Omega} \mathcal{L}\,^{s}_{\mathcal{B}}\, u(x) \ \mathcal{L}\,^{s}_{\mathcal{B}}\,v(x) \ dx.
\end{equation}
\end{remark}

\medskip
Now, the aim  is to characterize (via interpolation) the space $D(\mathcal{L}\,^{s}_{\mathcal{B}})$. 
To begin, we consider 
$u \in D\big(\,\mathcal{L}_{\mathcal{B}}\,\big)$, hence since $\mathcal{L}\,^{1/2}_{\mathcal{B}}$ is self-adjoint and from the definition of $\mathcal{L}_{\mathcal{B}}$ we have
$$
\begin{aligned}
\int_{\Omega} |\mathcal{L}\,^{1/2}_{\mathcal{B}} u(x)|^2 \ dx&=\int_{\Omega} \mathcal{L}\,^{1/2}_{\mathcal{B}}  u(x) \ \mathcal{L}\,^{1/2}_{\mathcal{B}}  u(x) \ dx=\int_{\Omega} \ \mathcal{L}_{\mathcal{B}}  u(x) \  u(x) \ dx
\\[5pt]
&=\int_{\Omega}-\dive\left( A(x)\nabla u(x) \right)u(x) \ dx=
\int_{\Omega}A(x)\nabla u(x) \cdot\nabla u(x) \ dx.
\end{aligned}
$$
On the other hand, using the uniform elliptic condition (see \eqref{unifellip}), we obtain
$$
  \Lambda_1 \int_{\Omega}|\nabla u(x)|^2 \ dx
  \leq \int_{\Omega} A(x)\nabla u(x) \cdot \nabla u(x) \ dx\leq \Lambda_2 \int_{\Omega}|\nabla u(x)|^2 \ dx, 
$$
where $\Lambda_2= \|A\|_\infty$. Therefore
\begin{equation}
 \Lambda_1 \Vert u \Vert_{H^1_{\Gamma_0}(\Omega)}^2\leq   \Vert 
 \mathcal{L}\,^{1/2 }_{\mathcal{B} } u \Vert_{L^2(\Omega)}^2\leq \Lambda_2 \Vert u \Vert_{H^1_{\Gamma_0}(\Omega)}^2,\label{normequi}
\end{equation}
which means the norm $\Vert\cdot\Vert_{1/2}$ is equivalent to the norm $\Vert \cdot\Vert_{H^1_{\Gamma_0}(\Omega)}$.  
Consequently, from the density of $D\big(\mathcal{L}_{\mathcal{B}}\big)$ in $D\big(\mathcal{L}\,^{1/2}_{\mathcal{B}}\big)$, and also in 
$H^1_{\Gamma_0}(\Omega)$, it follows that
$D\big(\mathcal{L}\,^{1/2}_{\mathcal{B}}\big)= H^1_{\Gamma_0}(\Omega)$. Similarly, we have the following 
\begin{proposition}
\label{represendoDominio}
Let $\Omega \subset \R^n$ be a bounded open set with Lipschitz boundary. 
If $s \in (0,1/2]$, then 
\begin{equation}
  D\big(\mathcal{L}\,^{s}_{\mathcal{B}}\big)= \left\{
 \begin{array}{rcl}
      &H^{2 s}(\Omega),&if\;\;\;\;\;0 < s < 1/4,
      \\[5pt]
      &\big[H^1_{\Gamma_0}(\Omega),L^2(\Omega)\big]_{1/2},&if\;\;\;\;\; s= 1/4,
      \\[5pt]
      &H^{2 s}_{\Gamma_0}(\Omega),&if\;\;\;\;\;1/4 <  s \leq 1/2.
\end{array}\right.\label{equalrepredominio}
\end{equation}
\label{eq:dominequiv}
\end{proposition}
\begin{proof}
The proof follows applying the discrete version of J-Method for interpolation, see \cite{MSireVazquez} and also \cite{Grubb}.
\end{proof}
%

Now, 
for each $s \in (0,1)$ we define conveniently the operators
$$
   \mathcal{K}_s:= \mathcal{L}^{-s}_{\mathcal{B}} \quad \text{and} \quad \mathcal{H}_s:= \mathcal{L}^{-s/2}_{\mathcal{B}}\equiv\mathcal{K}^{1/2}_s.
$$
Then we consider the following 
\begin{lemma}\label{lem:Kboundary}
Let $\;\Omega \subset \R^n$ be a bounded open set with Lipschitz boundary, $s\in(0,1)$ and $ u \in D(\,\mathcal{L}_\mathcal{B}\,)$, then
$\mathcal{K}_s u \in D(\,\mathcal{L}_\mathcal{B}\,)$. In particular, we have in trace sense
$$
\mathcal{K}_s  u= 0 \mbox{ on } \Gamma_0
\quad\mbox{ and }\quad
A\nabla \mathcal{K}_s u \cdot \nu =0 \mbox{ on } \Gamma_1.
$$
\end{lemma}
\begin{proof}
The proof follow directly from Proposition \ref{THMCARAC}, item (3).
\end{proof}

Here and subsequently, we denote $\mathbf{L}^2(\Omega)= \big(L^2(\Omega)\big)^n$.
Then we have the following important result. 
\begin{proposition}
\label{proKuHu}
Let $\;\Omega \subset \R^n$ be a bounded open set with Lipschitz boundary.
\begin{enumerate}
\item[$(1)$] 
If $u\in H^1_{\Gamma_0}(\Omega)$, then
$\nabla \mathcal{K}_s u \in \mathbf{L}^2(\Omega)$ and there exists 
$C_{\Omega}> 0$ such that
\begin{equation}
\label{nablaKunablau}
   \int_{\Omega}|\nabla\mathcal{K}_s u(x)|^2 \ dx
   \leq C_{\Omega}\int_{\Omega}|\nabla u(x)|^2 \ dx.
\end{equation}
Similarly, if $u\in H^1_{\Gamma_0}(\Omega)$, then
$\nabla \clg{H}_s u \in \mathbf{L}^2(\Omega)$ and
\begin{equation}
\label{nablaHunablau}
   \int_{\Omega}|\nabla\mathcal{H}_s u(x)|^2 \ dx
   \leq C_{\Omega}^{1/2} \int_{\Omega}|\nabla u(x)|^2 \ dx.
\end{equation}

\item[$(2)$] If $u\in H^1_{\Gamma_0}(\Omega)$, then
$$
\begin{aligned}
 &\Lambda_1\int_{\Omega}\vert \nabla \mathcal{H}_s u \vert^2 \ dx\leq   \int_{\Omega} A(x) \nabla \mathcal{K}_s u \cdot \nabla u \ dx 
 \leq \Lambda_2 \int_{\Omega}\vert \nabla \mathcal{H}_s u \vert^2 \ dx.
\end{aligned}
$$
\end{enumerate}
\end{proposition}

\begin{proof}
1. First, since $u\in H^1_{\Gamma_0}(\Omega)$, it is enough to consider $u \in D(\,\mathcal{L}_\mathcal{B}\,)$
and thus apply a standard density argument. 
To show item (1), from \eqref{normequi} we have
$$
\begin{aligned}
   \int_{\Omega}|\nabla\mathcal{K}_s u(x)|^2dx& \leq \Lambda_1^{-1}  \int_{\Omega}|\mathcal{L}_\mathcal{B}^{1/2}\mathcal{K}_s u(x)|^2dx
    = \Lambda_1^{-1} \sum_{k=1}^{\infty}\lambda_k\vert \langle \mathcal{K}_s u, \vp_k \rangle \vert^2
\\[5pt]
&=\Lambda_1^{-1} \sum_{k=1}^{\infty}\lambda_k\vert\lambda_k^{-s}\langle u, \vp_k \rangle\vert^2
\leq \Lambda_1^{-1}  \lambda_1^{-2s} \sum_{k=1}^{\infty}\lambda_k\vert\langle u, \vp_k \rangle\vert^2
\\[5pt]
   &= \Lambda_1^{-1} \lambda_1^{-2s} \!\! \int_{\Omega}|\mathcal{L}_\mathcal{B}^{1/2} u(x)|^2dx
\leq \Lambda_1^{-1} \Lambda_2 \ \lambda_1^{-2s} \!\! \int_{\Omega}|\nabla u(x)|^2dx< \infty,
\end{aligned}
$$ 
and analogously for $\nabla \clg{H}_A u$.

\medskip
2. Now, we prove item (2). Integrating by parts, we obtain
$$
\begin{aligned}
   \int_{\Omega}-\dive(A(x)\nabla\mathcal{K}_su(x))u(x)\,dx=&    \int_{\Omega}A(x)\nabla \mathcal{K}_su(x)\cdot\nabla u(x)\,dx
\\[5pt]   
   &-\int_{\Gamma}u(r)\,A(r)\nabla \mathcal{K}_su(r)\cdot\nu(r)\,dr. 
\end{aligned}
$$
We claim that, the boundary term is zero in the above equation. 
Let us recall that $\Gamma= \Gamma_0\cup\Gamma_1$, also since  $u\in D(\mathcal{L}_{\mathcal{B}})$, 
we have that $u=0$ on $\Gamma_0$. Moreover, from  Lemma \ref{lem:Kboundary} it
follows that $A \,\nabla \mathcal{K}_su\cdot\nu=0$ on $\Gamma_1$. Hence we conclude that the boundary term is zero. 
Therefore, we obtain
\begin{equation}
   \int_{\Omega}-\dive\big(A(x)\nabla\mathcal{K}_su(x)\big)\,u(x)\,dx=    \int_{\Omega}A(x)\nabla \mathcal{K}_su(x)\cdot\nabla u(x)\,dx\label{eq_pro_1}
\end{equation}

On the other hand, we observe that
 \begin{equation}
 \begin{aligned}
   \int_{\Omega}-\dive\big(A(x)\nabla\mathcal{K}_su(x)\big)u(x)\,dx&=    \int_{\Omega}\mathcal{L}_\mathcal{B} \left(\mathcal{K}_su(x)\right) \ u(x)dx
   \\[5pt]
   &=\int_{\Omega}\mathcal{L}_\mathcal{B}^{1-s}u(x) \, u(x)\,dx,
\end{aligned}\label{eq_pro_2}
\end{equation}
where we have used the definition of $\mathcal{L}_\mathcal{B}$ and $\mathcal{K}_s$. 
Then from \eqref{eq_pro_1}, \eqref{eq_pro_2} and since $\mathcal{L}_\mathcal{B}^{1-s}$ is self-adjoint (Proposition \ref{THMCARAC} ), it follows that
$$
\int_{\Omega}A(x)\nabla \mathcal{K}_su(x)\cdot\nabla u(x)\,dx= \int_{\Omega}|\mathcal{L}_\mathcal{B}^{(1-s)/2}u(x)|^2\,dx.
$$
Therefore, using the equivalence norm  \eqref{normequi} together with the definition of $\mathcal{H}_s u$, we have
$$
 \Lambda_1 \int_{\Omega}\vert \nabla \mathcal{H}_s u(x) \vert^2 \ dx
 \leq   \int_{\Omega}A(x)\nabla \mathcal{K}_s u(x)\cdot \nabla u(x)dx
  \leq \Lambda_2 \int_{\Omega}\vert \nabla \mathcal{H}_s u(x) \vert^2 \ dx.
$$
\end{proof}
\section{Initial Mixed-Boundary Value Problem} 
\label{SOLVAB}

The main issue of this section is to present the definition
of weak solutions for the 
initial mixed-boundary value problem \eqref{FTPME}, and 
then discuss in details in which sense
the initial mixed-boundary data will be considered, 
for any $s \in (0,1)$. 

\begin{definition}
\label{DEFSOL}
Given an initial data $u_0 \in L^\infty(\Omega)$ and $0< s <1$, a function 
$$
    u \in L^2\left((0,T); D\big(\,\mathcal{L}_{\mathcal{B}}^{(1-s)/2}\big)\right) \cap L^\infty(\Omega_T)
$$
is called a weak solution of  the initial mixed-boundary value problem \eqref{FTPME}, 
when $u$ satisfies
\begin{equation}
    \iint_{\Omega_T} u(t,x) \ (\partial_t\phi-A(x)\nabla \mathcal{K}_s u(t,x) \cdot \nabla \phi)dxdt
    +\int_{\Omega} u_0(x) \, \phi(0) dx= 0,\label{eq:thequiva}
\end{equation}
for each test function 
$\phi \in C^{\infty}_c \left([0,T) ; C^{\infty}_{\Gamma_0}(\bar{\Omega})\right)$. 
\end{definition}

One observes that, the above definition   makes sense. Indeed, the first and the last term in \eqref{eq:thequiva} is well defined, which is due to the fact that, 
$u$ and $u_0$ are bounded. The second term also works, it is enough to recall that $A(x)$ is bounded, and since for almost all $t \in (0,T)$,
$u(t) \in D\big(\,\mathcal{L}\,_{\mathcal{B}}^{(1-s)/2}\big)$, thus from item (3) in Proposition \ref{THMCARAC} 
and Proposition \ref{represendoDominio}, $ \mathcal{K}_s u(t) \in H^1_{\Gamma_0}(\Omega)$. Therefore, due to Proposition \ref{proKuHu} 
$$
\nabla \mathcal{K}_s u(t) \in \mathbf{L}^2(\Omega).
$$

\subsection{On the initial mixed-boundary data interpretation}
\label{BOCOMEN}

The aim of this section is to study the initial mixed-boundary datum interpretation, 
from the definition of weak solutions as presented by Definition \ref{DEFSOL}.
We start with the study of the mixed-boundary condition, and then the initial data
will be treat at the end of this section.

\medskip
To follow, we remark first that our definition of 
weak solutions is given for any $s \in (0,1)$, 
and hence it is not always possible to recover the boundary conditions in the trace sense.
Let us be more precise. The definition of a weak solution 
$$
    u\in L^2\left((0,T); D\big(\,\mathcal{L}\,_{\mathcal{B}}^{(1-s)/2}\big)\right) \cap L^\infty(\Omega_T)
$$
for \eqref{FTPME} is given by the integral equation 
\eqref{eq:thequiva}, where it is used a convenient space for the test functions, 
which give us some information about the mixed-boundary condition. 
Indeed, the homogeneous Dirichlet boundary condition is obtained by the space 
$D\big(\,\mathcal{L}\,_{\mathcal{B}}^{(1-s)/2}\big)$, and the Neumann boundary 
condition will be state via Coarea and Area Formulas. 

\medskip
Let $u$ be a solution of \eqref{FTPME} in the sense of Definition \ref{DEFSOL}. 
Firstly, we discuss  the Dirichlet condition, and it will be divided into three main steps:
 
1. If $0< s <1/2$ we have  
$$
u\in L^2\left((0,T); H^{1-s}_{\Gamma_0}(\Omega)\right),
$$
thanks to Proposition \ref{represendoDominio}. 
In particular, this space naturally encompass 
the Dirichlet boundary condition $u= 0$ on $\Gamma_0$, 
since the trace is well defined, see \eqref{trace}.

\medskip
2. Now, we consider $1/2<s<1$. In this case from Proposition \ref{represendoDominio}, we have
$$
u\in L^2\left((0,T); H^{1-s}(\Omega)\right).
$$
Here, the trace of $u$ on $\Gamma$ is not well defined, but we could give 
an interesting characterization. 
Indeed, applying Theorem 11.2 in \cite{LionsMagenes}, see p. 57, since 
for each $x \in \Omega$,
$\dist(x, \Gamma) \leq \dist(x, \Gamma_0)$, 
there exists a positive constant $C$, such that 
\begin{equation}
\label{CHARACT}
   \int_{\Omega}\frac{|u(\cdot,x)|^2}{(\dist(x, \Gamma_0))^{2(1-s)}} dx  \leq \frac{C}{2(1-s)}  \|u(\cdot)\|^2_{H^{1-s}(\Omega)}. 
\end{equation}
Now, since $\ \Gamma$ is a
$C^2-$boundary, there exists a sufficiently small $\delta> 0$ such that, each point 
$x \in \Omega_{\delta}:=\left\{ {x} \in \Omega :\ \dist(x, \Gamma)< \delta \right\}$ has a
unique projection ${r}=\mathbf{r}(x)$ on the
boundary $\Gamma$. Moreover, 
for every $x \in \Omega_\delta$ the Jacobian of the change of variables 
$$
   \Omega_\delta \ni x \leftrightarrow ({r},\tau) \in \Gamma \times (0,\delta)
   \quad \text{is equal to $\frac{D(x)}{D({r},\tau)}= 1 + O(\delta )$},
$$ 
where $\tau= \dist(x, \Gamma)$. Therefore, we obtain from \eqref{CHARACT}
\begin{equation}
\label{CHARACT10}
\begin{aligned}
  \int_0^\delta \!\! \int_{\Gamma_0}\frac{|u(\cdot,({r},\tau))|^2}{(\dist(({r},\tau), \Gamma_0))^{2(1-s)}} \ dr d\tau 
  &+  \int_0^\delta \!\! \int_{\Gamma_1}\frac{|u(\cdot,({r},\tau))|^2}{(\dist(({r},\tau), \Gamma_0))^{2(1-s)}} \ dr d\tau  
  \\[5pt] 
   &\leq \frac{C}{2(1-s)}  \|u(\cdot)\|^2_{H^{1-s}(\Omega)},
\end{aligned}
\end{equation}
and applying the Coarea Formula there exists a set of full mesure contained in $(0,\delta)$, such that, for each $\tau$ in this set  
$$
   \int_{\Gamma}\frac{|u(\cdot,({r},\tau))|^2}{(\dist(({r},\tau), \Gamma_0))^{2(1-s)}} \ dr
   \leq \frac{C}{2(1-s)}  \|u(\cdot)\|^2_{H^{1-s}(\Omega)}. 
$$
Moreover, for any $r \in \Gamma_0$ it follows that,
$\dist((r,\cdot), \Gamma_0)< \delta$. Hence we obtain from \eqref{CHARACT10} 
$$
   \lim \!\! \sup_{\delta \to 0^+} \Big(\delta^{2s-1} \ \frac{1}{\delta} \int_0^{\delta} \!\!\! \int_{\Gamma_0} |u(\cdot,({r},\tau))|^2 \,dr d\tau \Big) 
   \leq C,
$$
for some constant $C> 0$. Thus defining the following characterization
\begin{equation}
\label{CHARACTS}
    H^{1-s}_{\Gamma_{\!0(1-2s)}}(\Omega):=
    \left\lbrace f \in H^{1-s}(\Omega); \frac{1}{\tau} \int_0^{\tau} \!\!\! \int_{\Gamma_0} |f({r},\tau')|^2 \,dr d\tau'= O(\tau^{1-2s}) \right\rbrace,
\end{equation}
we have for almost all $t \in (0,T)$ that, $u(t) \in  H^{1-s}_{\Gamma_{\!0(1-2s)}}(\Omega)$, for any $1/2<s<1$. 

\medskip
3. The case $s=1/2$ is more delicate, since we do not have a precise indetification 
of the domain $D\big(\,\mathcal{L}\,_{\mathcal{B}}^{1/4}\big)$. 
Actually, from Proposition \ref{represendoDominio} and the second equation in \eqref{Eq_H_inclu}, we obtain
$$
H^{1/2}_{00}(\Omega)\subset D\big(\,\mathcal{L}\,_{\mathcal{B}}^{1/4}\big) \subset H^{1/2}(\Omega).
$$
First, we observe that the space $H^{1/2}(\Omega)$ 
does not have a well defined trace sense. 
On the other hand, there exists a notion of weak trace (see Theorem 11.7 in  \cite{LionsMagenes}) for 
$H^{1/2}_{00}(\Omega)$, but the spaces $H^{1/2}_{00}(\Omega)$ and $D\big(\,\mathcal{L}\,_{\mathcal{B}}^{1/4}\big)$ 
are not necessary equal. Although, we may follow the same strategy of item 2 above, and define the following 
characterization 
$$
    H^{1/2}_{\Gamma_{\!00}}(\Omega):=
    \left\lbrace f \in H^{1/2}(\Omega); \frac{1}{\tau} \int_0^{\tau} \!\!\! \int_{\Gamma_0} |f({r},\tau')|^2 \,dr d\tau'= O(1) \right\rbrace.
$$
Indeed, it is enough to observe that $D\big(\mathcal{L}_{\mathcal{B}}^{1/4}\big)$ is contained in $H^{1-s}(\Omega)$
for any $s \in [1/2, 1)$ and the right hand side of \eqref{CHARACT10} is uniformly bounded up to $s= 1/2$. 
Therefore for almost all $t \in (0,T)$, $u(t) \in  H^{1/2}_{\Gamma_{\!00}}(\Omega)$. 

\medskip
To finish the first part of this discussion, we study the Neumann boundary 
condition $\mathbf{q}(x,u)\cdot\nu=0$ on $\Gamma_1$, see \eqref{FTPME}, 
which is really complicated because it is composed by two terms, that is $u$ and $A \nabla\mathcal{K}_su\cdot \nu$.
In particular, we observe that $\mathbf{q}(x,u)\cdot\nu$ 
does not have trace on $\Gamma_1$ for any $0< s <1$. 
For instance, if $0<s<1/2$, it follows that $u(t)\in H_{\Gamma_0}^{1-s}(\Omega)\subset H^{1-s}(\Omega)$ a.e. in $(0,T)$, 
which implies that $u$ has trace on $\Gamma_1$ (not necessarily zero). Although,  there is no guarantee
that $A(x) \nabla \mathcal{K}_su\cdot \nu$ has trace on $\Gamma$, since $\mathcal{K}_su$ 
is not sufficiently regular. Similarly, if $1/2< s< 1$ then $\mathcal{K}_su$ 
is sufficiently regular to have trace on $\Gamma$, but 
as observed before we do not have trace for $u$. Thus 
the Neumann boundary condition is not well defined
in the strong sense in any case.

\medskip
On the other hand, Definition \ref{DEFSOL} is sufficiently robust to give a sense of the Neumann boundary condition on $\Gamma_1$. 
More precisely, we state this boundary condition in a weak sense, written as limits of integrals on $(0,T)\times\Gamma_1$. Indeed, we prove that any solution $u$ in the sense of Definition \ref{DEFSOL}, satisfies 
\begin{equation}
\ess\lim_{\tau \to 0^+}\int_0^T\int_{\Gamma_1} \mathbf{q}\left(\Psi_\tau(r), u(t,\Psi_\tau(r)) \right)\cdot \nu_\tau(\Psi_\tau(r)) \ \phi(t,r) \,dr\,dt= 0,\nonumber
\end{equation}
where $\Psi_\tau(r):=r-\tau\nu(r)$, and $\nu_\tau$ is the unit outward normal field on $\Psi_\tau(\Gamma)$,
see Appendix. 

\medskip
To prove the above sentence, we consider the following sets: Let $\clg{F}$ be a countable dense subset of $C^{\infty}_c\left((0,T); C^1_{\Gamma_0}(\bar{\Omega})\right)$.
For each $\gamma \in \clg{F}$, we define the set of full measure in $(0,1)$ by
$$
 \text{$F_\gamma= \Big \{ \tau\in(0,1) / \tau$ is a Lebesgue point of $\mathbf{J}(\tau) \Big \}$}, 
$$
where $\mathbf{J}(\tau)$ is given by
$$
\int_0^{T}\!\!\int_{\Gamma_1}\mathbf{q}\left(\Psi_\tau(r), u(t,\Psi_\tau(r)) \right)\cdot\nu_\tau(\Psi_\tau(r))J[\Psi_\tau(r)]\gamma(t,r)\,dr \ dt,
$$
where $J[\Psi_\tau]$ is the Jacobian of $\Psi_\tau$. Then, we consider 
$$
F:=\bigcap_{\gamma\,\in\, \clg{F}}F_\gamma,
$$
which is also a set of full measure in $(0,1)$. 
 
\begin{proposition}[Neumann condition]
\label{pro:NC} Let $u$ be a weak solution
 for the 
initial mixed-boundary value problem \eqref{FTPME}, 
in the sense of Definition \ref{DEFSOL}. Then, for each $\gamma \in H^1_0(0,T; H^1_{\Gamma_0}(\Omega))$
\begin{equation}
\ess\!\!\lim_{\tau \to 0^+}\int_0^{T}\int_{\Gamma_1}\!\!\!\mathbf{q}\left(\Psi_\tau(r), u(t,\Psi_\tau(r)) \right)\cdot\nu_\tau(\Psi_\tau(r))\gamma(t,r)\,dr\,dt= 0,\nonumber
\end{equation}
where $\Psi_\tau(r):=r-\tau\nu(r)$
and $\nu_\tau$ is the unit outward normal field in $\Psi_\tau(\Gamma)$.
\end{proposition}
\begin{proof} 
First, we define $S:= \Psi(F \times \Gamma)$ and consider
$$
 \phi(t,x)=\left\lbrace
 \begin{aligned}
 \gamma(t,\Psi_{h(x)}^{-1}(x))\zeta_j(h(x)), &\quad \text{for $x \in S$}, 
 \\[5pt]
0\;, & \quad \text{for $x \in \Omega \setminus S$},
 \end{aligned}\right.
$$
where $\gamma \in \clg{F}$, $\zeta_j(\tau)=H_j(\tau+\tau_0)-H_j(\tau-\tau_0)$, with $\tau_0\in F$.
Therefore, from \eqref{eq:thequiva} with $\phi(t,x)$ as test function, and applying the Coarea Formula for the function $h$, we have
$$
 \begin{aligned}
 &\int^1_0\zeta_j(\tau)\int^T_0\int_{\Psi_\tau(\Gamma)}u(t,r)\partial_t\gamma(t,\Psi_{\tau}^{-1}(r))\,d\clg{H}^{n-1}(r)dtd\tau
 \\[5pt]
  &=\int^1_0\zeta_j(\tau)\int^T_0\int_{\Psi_\tau(\Gamma)}\mathbf{q}(r, u(t,r))
  \cdot \nabla\gamma(t,\Psi_{h(x)}^{-1}(x))(r)\,d\clg{H}^{n-1}(r)dtd\tau
 \\[5pt]
 &+\int^1_0\zeta_j'(\tau)\int^T_0\int_{\Psi_\tau(\Gamma)}\mathbf{q}(r, u(t,r))
 \cdot \nu_\tau(r) \gamma(t,\Psi_{\tau}^{-1}(r))\,d\clg{H}^{n-1}(r)dtd\tau,
 \end{aligned}
$$
where we have used \eqref{eq:nablah} and $\nabla h$ is parallel to $\nu_\tau$ $\mathcal{H}^{n-1}$ a.e on $\Psi_\tau(\Gamma)$.

\medskip
Then, using the Area formula for the function $\Psi_\tau$ and  passing to the limit in the above equation as 
$j\to\infty$, recall that $\tau_0$ is a Lebesque point of $\mathbf{J}(\tau)$, moreover $\zeta_j(t)$ 
converges pointwise to the characteristic function of the interval $[-\tau_0,\tau_0)$ and $\gamma(t,\cdot)=0$ on $\Gamma_0$, we obtain
\begin{equation}
  \mathbf{J}(\tau_0)=\int_0^{\tau_0}\Phi(\tau)d\tau,\label{eq:J(tau)}
\end{equation}
 for all $\tau_0\in F$ and $\gamma\in\clg{F}$, where $\Phi(\tau)$ is given by
 $$
 \int^T_0 \!\! \int_{\Psi_\tau(\Gamma)}u(t,r) \big(\partial_t \gamma(t,\Psi_{\tau}^{-1}(r))
 -A(r)\nabla\mathcal{K}_su(t,r)\cdot\nabla\gamma(t,\Psi_{h(\cdot)}^{-1}(\cdot))(r)\big) d\clg{H}^{n-1}(r)dt.
 $$
 
 \smallskip
On the other hand, since $\clg{F}$ is dense in $C^{\infty}_c\big((0,T); C^1_{\Gamma_0}(\bar{\Omega})\big)$, we have that \eqref{eq:J(tau)}
holds for $\gamma\in C^{\infty}_c\big((0,T); C^1_{\Gamma_0}(\bar{\Omega})\big)$. Then, for each $\tau\in F$ we have
$$
\left| \mathbf{J}(\tau)\right|\leq C\left| \Psi\left(\,(0,\tau)\times\Gamma\,\right)\right|,
$$
where $C$ is a positive constant, which does not depend on $\tau$. 
Moreover, we know that $J[\Psi_\tau] \to 1$ as $\tau\to 0^+$.
Therefore, applying the Dominated Convergent Theorem we obtain  
\begin{equation}
\ess\!\!\lim_{\tau \to 0^+}\!\!\int_0^{T}\!\!\!\int_{\Gamma_1}\!\!\!
\mathbf{q}(\Psi_\tau(r), u(t,\Psi_\tau(r)))\cdot\nu_\tau(\Psi_\tau(r))\gamma(t,r) \,drdt= 0,\nonumber
\end{equation}
which completes the proof.
\end{proof}

\bigskip
To finish this section, we characterize the initial boundary condition from Definition  \ref{DEFSOL}.
For this purpose, let $\clg{E}$ be a countable dense subset of $C^{1}_{\Gamma_0}(\bar{\Omega})$.
For each $\zeta \in \clg{E}$, we define the set of full measure in $(0,T)$ by
$$
 \text{$E_\zeta:= \Big\{ t \in(0,T) / \, t $ is a Lebesgue point of $ \mathrm{I}(t)= \int_{\Omega}u(t,x)\zeta(x)dx \Big\}$},
$$
and consider
$$
E:=\bigcap_{\zeta\,\in\, \clg{E}}E_\zeta,
$$
which is a set of full measure in $(0,T)$.
\begin{proposition}[Initial condition]
 Let $u$ be a weak solution
 for the 
initial mixed-boundary value problem \eqref{FTPME}, 
in the sense of Definition \ref{DEFSOL}. Then for all $\zeta\in L^1(\Omega)$
\begin{equation}
   \ess \lim_{t \to 0^+} \int_\Omega u(t,x) \zeta(x)\, dx=  \int_\Omega  u_0(x) \zeta(x)\, dx.\label{eq:initalcondition}
\end{equation}
\end{proposition}
\begin{proof}
We give only the main ideas of the proof (for more details see \cite{Hua}). Let us consider $\phi(t,x)= \gamma_j(t)\zeta(x)$, 
$\gamma_j(t)=H_j(t+t_0)-H_j(t-t_0)$ for any $t_0 \in E$ (fixed), and 
$\zeta \in \clg{E}$. 
Then,  substituting $\phi$ into \eqref{eq:thequiva} and
passing to the limit as $j \to \infty$, ( $t_0$ is Lebesque point of 
$\mathrm{I}(t)$ ), we obtain
\begin{equation}
   \mathrm{I}(t_0) =\int_{\Omega} u_0(x)\zeta(x)dx -\int^{t_0}_0\int_{\Omega}u(x)A(x)\nabla\mathcal{K}_su(x)\cdot\nabla\zeta(x)dxdt,\label{INT100}
\end{equation} 
where we have used the Dominated Convergence Theorem. Since $t_0 \in E$ is arbitrary, and in view of the density of $\clg{E}$ 
in $L^1(\Omega)$, the proof
follows. 
\end{proof}

\section{Main Result}
\label{MainResult}
The main result of this section is to show a weak solution of \eqref{FTPME}. To this end, we have the following 

\begin{theorem} [Main Theorem]
\label{Thprincipal}
Let $u_0 \in L^{\infty}(\Omega)$
be a non-negative function. 
Then, there exists a weak solution $u \in L^2\big((0,T);D\big(\,\mathcal{L}\,_{\mathcal{B}}^{(1-s)/2}\big)\big)\cap L^{\infty}(\Omega_T)$ 
of the initial mixed-boundary value problem \eqref{FTPME}.
\end{theorem}
The proof of this result is given in the next sections. 

\subsection{Anisotropic parabolic approximation}
\label{Parabolic regularization}
In this subsection, we introduce and  study
the approximate parabolic problem with  $\delta, \mu \in (0, 1)$, given by
\begin{eqnarray}
     \partial_t u_{\mu,\delta}
     -\delta \ \dive(A(x)\nabla u_{\mu,\delta})&=&
     \dive(\mathbf{q}_{\mu}(x,u_{\mu,\delta})) 
     \;\hspace{0,8cm}\;\mbox{in}\;\Omega_T,\label{eq:aproxequation1}
     \\[5pt]
{u_{\mu,\delta}} &=&u_{0\delta} \hspace{2.5cm}\mbox{in}\; \{t=0\} \times \Omega,\label{eq:aproxequation2}\\[5pt]
 u_{\mu,\delta}&=& 0 \hspace{2,8cm}\mbox{on}\; (0,T) \times \Gamma_0,\label{eq:aproxequation3}
 \\[5pt]
 \delta A\nabla u_{\mu,\delta} \cdot \nu&=& -\mathbf{q}_{\mu}(x,u_{\mu,\delta}) \cdot \nu \hspace{0,6cm}\mbox{on}\; (0,T) \times \Gamma_1,\label{eq:aproxequation4}
\end{eqnarray}
where $\mathbf{q}_{\mu}(x,u):= (\mu+u) A(x) \nabla \mathcal{K}_s u$, and $u_{0\delta}$ is a non-negative 
regularized initial data 
such that 
\begin{equation}
\text{$u_{0,\delta} \to u_0 $ strongly in $L^{1}(\Omega)$
as $\delta \to 0$, \quad
$\|u_{0,\delta}\|_{L^\infty} \leq \|u_{0}\|_{L^\infty}$,}\nonumber
\end{equation}
and satisfying suitable compatibility
conditions.

\medskip
Now,  we make use of the well known results of existence, uniqueness and uniform $L^{\infty}$ bounds 
for parabolic problems with mixed boundary conditions. Therefore, applying Theorem \ref{EXUNIQAPP} in Appendix,
for each $\delta,\;\mu> 0$, 
there exists a unique, namely here strong solution, 
$$
\begin{aligned}
u_{\mu,\delta}&\in C([0,T); H^1_{\Gamma_0}(\Omega))\cap L^2((0,T); H^2(\Omega^{\prime}))  \cap L^\infty(\Omega_T), 
\\[5pt]
\partial_t u_{\mu,\delta} &\in L^2\left(\Omega_{T}\right),
\end{aligned}
$$
for each $\Omega'$ compactly contained in $\Omega$. 
Moreover, one observes that conditions \eqref{eq:aproxequation3} and \eqref{eq:aproxequation4} are satisfied 
in the sense of trace. 

\medskip
The following theorem investigates the properties of the solution
$u_{\mu,\delta}$ to the (anisotropic) parabolic perturbation 
\eqref{eq:aproxequation1}--\eqref{eq:aproxequation4} for fixed $\delta,\mu \in (0, 1)$. 
\begin{theorem}
\label{Theaprox}
For each $\mu,\delta>0$, let 
$u= u_{\mu,\delta}$
be the unique strong solution of  
\eqref{eq:aproxequation1}--\eqref{eq:aproxequation4}. Then, $u$ satisfies:
\begin{enumerate}
\item[$(1)$] For all $\phi\in C_c^{\infty}([0,T):C^{\infty}_{\Gamma_0}(\bar{\Omega}))$,
\begin{equation}
\label{eq:thequiaproxsol}
\hspace{-0,57cm}\begin{aligned}
    \iint_{\Omega_T} (u(t,x)\partial_t \phi(t,&x)-\delta A(x)\nabla u\cdot\nabla\phi(t,x)
    ) \ dxdt 
    + \int_{\Omega}u_{0\delta}(x) \ \phi(0,x) \ dx
 \\[5pt]
    &=\iint_{\Omega_T} (\mu+u(t,x)) A(x) \nabla \mathcal{K}_su(t,x) \cdot \nabla\phi(t,x) \ dxdt.
\end{aligned}
\end{equation}

\item[$(2)$] For all $(t,x)\in\Omega_T$, we have 
\begin{equation}
\label{eq:+ebounded}
   0\leq  u(t,x)+\mu\leq\Vert u_0\Vert_{L^\infty},
\end{equation}  
and the conservation of the ``total mass"
\begin{equation}
\label{eq:consermassaprox}
\int_{\Omega}u(t,x) \ dx = \int_{\Omega} u_{0 \delta}(x) \ dx
\leq \|u_0\|_{L^\infty} |{\Omega}|.
\end{equation}
\end{enumerate}
\end{theorem}
\smallskip
\begin{proof}
1. Let us show \eqref{eq:thequiaproxsol}. First, we observe that, the equation \eqref{eq:aproxequation1} is verified 
for almost all points $(t,x) \in (0,T) \times \Omega^\prime$, for each $\Omega^\prime$ compactly contained in $\Omega$. 
Therefore, we multiply \eqref{eq:aproxequation1} by $\phi(t,x) \, (1-\zeta_j(h(x)))$ and integrate in $\Omega_T$,  
where $\phi\in C_c^{\infty}([0,T); C^{\infty}_{\Gamma_0}(\bar{\Omega}))$, and $\zeta_j(h(x))$ is taken as in the 
proof of Proposition \ref{pro:NC}. We are not going to reproduce here all the details given at Section \ref{BOCOMEN},
and from now on we omit this procedure. One remarks that, the support of $(1-\zeta_j(h(x))) \subset \Omega$.   
Then, after integration by parts we obtain 
\begin{equation}
\label{EQAPP}
\begin{aligned}
 & \int_0^T  \!\!\! \int_{\Omega} \big\{  -u \partial_t\phi+\delta A(x)\nabla u\cdot \nabla\phi
   + (\mu+u)A(x)\nabla\mathcal{K}_s u \cdot \nabla \phi \big\} (1-\zeta_j) \ dx dt
\\[5pt] 
&\quad =  \int_{\Omega} u_{0\delta} \ \phi(0) \, (1-\zeta_j) \ dx
+ \int_0^T \!\!\! \int_{\Gamma} \phi \, (1-\zeta_j) \,\left( \delta  A(r)\nabla u +\mathbf{q}_{\mu}(r,u) \right)\cdot \nu\, dr dt
 \\[5pt]
 &\quad + \int^1_0 (-\zeta_j'(\tau)) \int^T_0 \!\!\! \int_{\Psi_\tau(\Gamma)} \phi \left( \delta  A(r)\nabla u +\mathbf{q}_\mu(r, u) \right)
 \cdot \nu_\tau(r) \,dr dt d\tau,
\end{aligned}\nonumber
\end{equation}
where we have used the Coarea Formula for the function $h$ in the third integral in the right hand side of the above equation.
Thus, applying the Area formula for the function $\Psi_\tau$, passing to the limit as $j \to \infty$ and making $\tau_0 \to 0^+$, 
 we have
$$
\begin{aligned}
 \int_0^T  \!\!\! \int_{\Omega} &\big\{  -u \partial_t\phi+\delta A(x)\nabla u\cdot \nabla\phi
   + (\mu+u)A(x)\nabla\mathcal{K}_s u \cdot \nabla \phi \big\}\ dxdt
\\[5pt] 
&=  \int_{\Omega} u_{0\delta} \ \phi(0) \ dx 
+ 2 \int_0^T \!\!\! \int_{\Gamma} \phi \,\left( \delta  A(r)\nabla u +\mathbf{q}_{\mu}(r,u) \right)\cdot \nu\, dr dt. 
\end{aligned}
$$
Finally, we stress that the boundary term 
$$
   \int_0^T \!\!\! \int_{\Gamma} \phi \,\left( \delta  A(r)\nabla u +\mathbf{q}_{\mu}(r,u) \right)\cdot \nu\, dr dt= 0.
$$
Indeed, due to $\Gamma=\Gamma_0\cup\Gamma_1$, $\phi= 0$ on 
$\Gamma_0$, and $\left( \delta  A(r)\nabla u +\mathbf{q}_{\mu}(r,u) \right)\cdot \nu=0$ on $\Gamma_1$,
see \eqref{DEFTESTFUNC} and \eqref{eq:aproxequation4} respectively. 

\medskip
2. To show the assertion \eqref{eq:+ebounded}, we multiply \eqref{eq:aproxequation1} by $\varphi_{\varepsilon}^{\prime}(u)$ 
and integrate in $\Omega_t= (0,t)\times\Omega$, $0< t \leq T$, where
\[
\varphi_{\varepsilon}(z)=\left\{\begin{array}{ll}
\left((z+\mu)^{2}+\varepsilon^{2}\right)^{1 / 2}-\varepsilon, & \text { for } z \leq -\mu,
\\[5pt]
0, & \text { for } z \geq -\mu,
\end{array}\right. 
\]
which converges to $|z+\mu|^{-}:=\min\{z+\mu,0\}$ as $\varepsilon\to0^+$. Hence from the properties of $\varphi_\varepsilon$, we obtain
\[
\begin{aligned}
\int_{\Omega} \varphi_{\varepsilon}(u(t)) d x+\iint_{\Omega_t}\varphi_{\varepsilon}^{\prime \prime}(u)(\mu+u(x))A(x)\nabla\mathcal{K}_s u\cdot \nabla u \  d x d \tau 
\\[5pt]
+\delta \iint_{\Omega_t} \varphi_{\delta}^{\prime \prime}(u_\varepsilon)\,A(x)\nabla u\cdot\nabla u\, d x d \tau= 0,
\end{aligned}
\]
where we have used that, $u_0 \geq 0$, the boundary conditions in \eqref{eq:aproxequation3}-\eqref{eq:aproxequation4}, 
and $\varphi^{\prime}_\varepsilon(0)=0$. On the other hand, we observe
\[
\begin{aligned}
\varphi_{\varepsilon}^{\prime \prime}(u)(\mu+u(x))A(x)\nabla\mathcal{K}_s u\cdot \nabla u&+\delta A(x)\nabla u\cdot\nabla u\, \varphi_{\varepsilon}^{\prime \prime}(u)
 \\[5pt]
& \geq\left\{-|\mu+u(x)||A(x)\nabla\mathcal{K}_s u||\nabla u|+\delta\Lambda_1|\nabla u|^{2}\right\} \varphi_{\varepsilon}^{\prime \prime}(u) 
\\[5pt]
& \geq-\frac{1}{4 \delta\Lambda_1}(\mu+u)^{2}|A(x)\nabla \mathcal{K}_s u|^2 \varphi_{\varepsilon}^{\prime \prime}(u)
\\[5pt]
& \geq-\frac{\varepsilon}{4 \delta\Lambda_1}|A(x)\nabla\mathcal{K}_s u|^2,
\end{aligned}
\]
where we have used the uniform ellipticity and $(u+\mu)^2\varphi_{\varepsilon}^{\prime \prime}(u)\leq\varepsilon$. Consequently,
\[
\int_{\Omega} \varphi_{\varepsilon}(u(t)) d x\leq\frac{\varepsilon}{4 \delta\Lambda_1}\int_{\Omega_t}|A(x)\nabla\mathcal{K}_s u(\tau,x)|^2\,dx\,d\tau.
\]
Then passing the limit as $\varepsilon\to0^+$, we get
$$
\int_\Omega|u(t,x)+\mu|^{-} \ dx\leq 0,
$$
thus $|u(t,x)+\mu|^{-}=0$. Similarly, we can show that $|u(t,x)+\mu-\|u\|_{\infty}|^{+}=0$, therefore \eqref{eq:+ebounded}
is proved.

\medskip
3. It remains to prove \eqref{eq:consermassaprox}. We multiply \eqref{eq:aproxequation1} by $\xi_k(x)$ (see Appendix), and integrate
over $\Omega$. Then, after integration by parts and due to $\xi_k=0$ on $\Gamma$, we have
$$
\begin{aligned}
    \dfrac{\partial }{\partial t}\int_\Omega u(t,x) \xi_k(x) \, dx
 =&-\int_{\Omega}\delta \ A(x)\nabla u(t,x)\cdot\nabla \xi_k(x) \ dx
\\[5pt]
 & -\int_{\Omega} (\mu+u(t,x))A(x)\nabla\mathcal{K}_s u(t,x)  \cdot\nabla \xi_k(x) \ dx.
\end{aligned}
$$
Now, we integrate the above equation over $(0, t)$
\begin{equation}
\begin{aligned}
\int_{\Omega} \big(u(t,x)&-u_{0,\delta}(x) \big) \xi_k(x) \ dx =-\int_0^t\int_{\Omega}\delta \ A(x)\nabla u(t,x)\cdot\nabla \xi_k(x) \ dx
\\[5pt]
&-\int_0^t\int_{\Omega}(\mu+u(t,x))A(x)\nabla\mathcal{K}_s u(t,x))\cdot\nabla \xi_k(x) \ dxdt'
\\[5pt]
&=-I_1-I_2, 
\end{aligned}
\end{equation}
with the obvious notation. Let us observe the $I_2$ term, we have
$$
   | I_2| \leq (\|u\|_\infty + 1) \|A\|_\infty  
   \Big(\iint_{\Omega_T} |\nabla\mathcal{K}_s u(t,x)|^2 dxdt \Big)^{1/2}
   \Big(\iint_{\Omega_T} |\nabla \xi_k(x)|^2 dxdt \Big)^{1/2},
$$
where we have used H\"{o}lder's inequality and the uniform estimates for $u(t,x)$, $A(x)$. 
Therefore, applying Lemma \ref{Lemma:xik} we obtain
$$
    \lim_{k \to \infty} I_2= 0. 
$$
Similarly, we have that $I_1$ goes to zero as $k \to \infty$. 
Then, passing to the limit as $k\to \infty$ in \eqref{EQAPP}, and again 
applying Lemma \ref{Lemma:xik} we get \eqref{eq:consermassaprox}. 
Hence the proof of the Theorem \ref{Theaprox} is complete.
\end{proof}

Now, let us consider two important estimates of the solution $u_{\delta,\mu}$ for 
the initial mixed-boundary valued problem \eqref{eq:aproxequation1}--\eqref{eq:aproxequation4}, with fixed $\delta, \mu\in (0, 1)$.

\begin{proposition}[First energy estimate]\label{eq:Fu}
Let $u= u_{\mu,\delta}$ be the unique strong solution of  
\eqref{eq:aproxequation1}--\eqref{eq:aproxequation4}. Then, 
 for all $t\in (0,T)$,
\begin{equation}
\begin{aligned}
\int_{\Omega}\eta(u(t))\,dx&+\Lambda_1\delta\int_0^t\int_{\Omega} \ \dfrac{\vert \nabla u \vert^2}{\mu+u} \ \,dx\,dt\,
\\[5pt]
&+\Lambda_1\int_0^t\int_{\Omega}|\nabla \mathcal{H}_s u|^2\,dxdt\leq   \int_{\Omega}\eta(u_{0\delta})dx,
\end{aligned}\label{eq:firstenergy}
\end{equation}
where $\eta(\lambda):= (\lambda+\mu) \log (1+(\lambda/\mu))-\lambda$, $(\lambda\geq 0)$.
\end{proposition}
\begin{proof}
First, we multiply \eqref{eq:aproxequation1}
by $\eta'(u)$ and integrate on $\Omega$. Then, after integration by parts, we have
$$
\begin{aligned}
   \dfrac{\partial}{\partial t}\int_{\Omega}\eta(u)dx=&-\delta\int_{\Omega}
   \dfrac{1}{\mu+u} \  A(x)\nabla u\cdot\nabla u  \ dx
   - \int_{\Omega}A(x)\nabla \mathcal{K}_s u\cdot\nabla u\,dx
   \\[5pt]
   &+\int_{\Gamma}\eta'(u(r)) \ \left(\delta A(r)\nabla u(r)+ \mathbf{q}_{\mu}(r,u)\right)\cdot\nu \ dr.
\end{aligned}
$$
One observes that, the boundary terms are zero. Indeed, the proof is similar to Theorem \ref{Theaprox},
where the important point here is that $\eta'(0)=0$ and $u=0$ on $\Gamma_0$. Therefore, the boundary terms are zero.
Then, we integrate over $(0,t)$, for all $0 < t < T$, to obtain 
$$
\begin{aligned}
   \int_{\Omega}\eta(u(t))dx &+ \delta \int_0^t\int_{\Omega}
   \dfrac{1}{\mu+u(t,x)} \ A(x) \nabla u(t,x) \cdot \nabla u(t,x)  \ dx dt
   \\[5pt]
   &+ \int_0^t\int_{\Omega}A(x) \nabla \mathcal{K}_s u(t,x) \cdot \nabla u(t,x) \,dx dt= \int_{\Omega}\eta(u_0) \ dx.
\end{aligned}
$$

On the other hand, due to the uniform ellipticity condition, we have 
$$
\Lambda_1\int_0^t\int_{\Omega} \ \dfrac{|\nabla u(t,x) |^2}{\mu+u(t,x)}  \ dx dt 
\leq\int_0^t \int_{\Omega} \dfrac{1}{\mu+u(t,x)} \  A(x)\nabla u(t,x)\cdot \nabla u(t,x)  \ dx dt. 
$$
For the third term in the left hand side, we use Proposition \ref{proKuHu} ($u\in H^1_{\Gamma_0}(\Omega)$), 
which establishes the first energy estimate.
\end{proof}

As a consequence of this last result, we obtain
\begin{corollary}\label{cor:gradlimi} Under the assumptions of the Proposition \ref{eq:Fu}, we have that $u=u_{\delta,\mu}$ satisfies
\begin{equation}
\begin{aligned}
\delta\|\nabla u\|^2_{L^2(\Omega_T)} &\leq \|u_0\|_{\infty} \ \eta(\|u_0\|_{\infty}) \ |\Omega| \ \Lambda_1^{-1},\quad\mbox{ and }
\\[5pt]
\|\nabla \mathcal{H}_s u\|^2_{L^2(\Omega_T)} &\leq \eta(\|u_0\|_{\infty}) \ |\Omega| \ \Lambda_1^{-1}, 
\end{aligned}\label{eq:gradlimi}
\end{equation}
where $|\Omega|$ is the Lebesgue measure of the set $\Omega$.
\end{corollary}
\begin{proof}
We only provide the proof for the first inequality in \eqref{eq:gradlimi}, the other one is similar. From \eqref{eq:firstenergy} we have
$$
\frac{\Lambda_1 \delta}{ \|u_0\|_{\infty}} \int_0^t\int_{\Omega} \ |\nabla u(t,x)|^2 \,dx dt
\leq   \int_{\Omega}\eta(u_{0\delta}(x)) dx,
$$
where we have used \eqref{eq:+ebounded}. 
Moreover, since $\eta'(\lambda)> 0$, $(\lambda \geq 0)$, if follows that 
$\eta(\lambda)$ is an increasing function, hence $\eta(u_{0\delta}(x)) \leq \eta(\|u_0\|_{\infty})$
for almost all $x \in \Omega$. 
Consequently, we obtain 
$$
\int_{\Omega}\eta(u_{0\delta}(x))dx \leq \eta(\|u_0\|_{\infty})|\Omega|, 
$$
which completes the proof.
\end{proof}
\begin{proposition}[Second energy estimate]\label{eq:2energystimate} Under the conditions stated above, we have that $u= u_{\mu,\delta}$ satisfies
\begin{equation}
\begin{aligned} 
\frac{1}{2}\int_{\Omega}&|\mathcal{H}_s u(t_2,x)|^2\,dx+\Lambda_1\delta \int_{t_1}^{t_2}\int_{\Omega}|\nabla\mathcal{H}_s u|^2\,dxdt
 \\[5pt]
&+\Lambda_1\int_{t_1}^{t_2}\int_{\Omega} (\mu+u) \left|\nabla \mathcal{K}_s u\right|^2\,dxdt
  \leq \frac{1}{2}\int_{\Omega}|\mathcal{H}_s u(t_1,x)|^2 \ dx,
\end{aligned}\label{eq:secondenergy}
\end{equation}
 for all  $0 \leq t_1 < t_2 < T$.
\end{proposition}
\begin{proof}
First, we multiply \eqref{eq:aproxequation1} by $\mathcal{K}_su$, and
integrate in $\Omega$. Then, we have
$$
\begin{aligned}
\int_{\Omega}\frac{\partial u}{\partial t}\mathcal{K}_su \ dx
=&-\delta\int_{\Omega}A(x)\nabla u\cdot\nabla\mathcal{K}_su \ dx 
-\int_{\Omega} (\mu+u)A(x)\nabla \mathcal{K}_s u\cdot\nabla \mathcal{K}_su\,dx
\\[5pt]
&+\int_{\Gamma}\mathcal{K}_su \ \left(\delta A(r)\nabla u + \mathbf{q}_{\mu}(r,u)\right)\cdot\nu \ dr.
\end{aligned}
$$
One observes that,
$u(t) \in H^1_{\Gamma_0}(\Omega)$ for each $t \in [0,T)$, thus 
by Proposition \ref{THMCARAC} it follows that $\mathcal{K}_s u(t)= 0$ on $\Gamma_0$. Hence, from the same ideas used above, we have that the boundary terms are zero. 
Then, integrating over $0 \leq t_1 < t_2 < T$, we obtain
$$
\begin{aligned}
\frac{1}{2}\int_{\Omega}|\mathcal{H}_su(t_2,x)|^2dx+\delta\int_{t_1}^{t_2}\int_{\Omega}A(x)\nabla u\cdot\nabla\mathcal{K}_su\,dx dt
&
\\[5pt]
+ \int_{t_1}^{t_2}\int_{\Omega} (\mu+u)A(x)\nabla \mathcal{K}_s u\cdot\nabla \mathcal{K}_su\,dx dt&=
\frac{1}{2}\int_{\Omega}|\mathcal{H}_su(t_1,x)|^2dx.
\end{aligned}
$$
From the uniform ellipticity condition, we have and estimate for the third term of the left hand side 
$$
\Lambda_1 
\int_{t_1}^{t_2}\int_{\Omega}(\mu+u)|\nabla \mathcal{K}_su |^2  \ dx\leq\int_{t_1}^{t_2}\int_{\Omega}(\mu+u) A(x)\nabla \mathcal{K}_su\cdot\nabla \mathcal{K}_s u  \ dx
$$
and for the second term, we use  Proposition \ref{proKuHu} ($u\in H^1_{\Gamma_0}(\Omega)$).
Therefore we get the second energy estimate \eqref{eq:secondenergy}.
\end{proof}

Finally, we consider the following 
\begin{proposition}
\label{eq:L2Hs-1}
Under the above conditions, we have 
 for all $v\in H^1_{\Gamma_0}(\Omega)$
\begin{equation}
\int^T_0 \langle \partial_t u(t), v\rangle dt=-\delta\iint_{\Omega_T}A(x)\nabla u\cdot\nabla v\,dxdt 
+ \iint_{\Omega_T} (\mu+u)A(x)\nabla \mathcal{K}_su\cdot \nabla v \,dxdt.\label{eq:3estimativa}
\end{equation}
where $\langle\cdot,\cdot\rangle$ denotes the pairing between $\left(H^1_{\Gamma_0}(\Omega)\right)^*$ and $H^1_{\Gamma_0}(\Omega)$.
\end{proposition}
\begin{proof}
The proof follows applying the same techniques considered before, so it is omitted. 
\end{proof}

\subsection{Proof of Main Theorem}

Here we pass to the limit in \eqref{eq:thequiaproxsol}, as the two parameters $\delta$, $\mu$ go to zero. 
To this end, we use the first and the second energy estimates together with the Aubin-Lions' Theorem. 

\subsubsection{Limit transition $\delta\to 0^+$}

As a first step, we define $u_\delta:= u_{\mu,\delta}$ (fixing $\mu> 0$). 
The main result in this section is the following
\begin{proposition}
\label{Thprinc}
Let $\{u_\delta\}_{\delta>0}$ be the strong solutions of \eqref{eq:aproxequation1}--\eqref{eq:aproxequation3}. 
Then, there exists a subsequence of $\{u_\delta\}_{\delta>0}$, which 
weakly converges to some function $u\in  L^2\big((0,T); D\big(\mathcal{L}_{\mathcal{B}}^{(1-s)/2}\big)\big) \cap L^{\infty}(\Omega_T)$,
satisfying
\begin{equation}
\begin{aligned}
\iint_{\Omega_T}u(t,x)\partial_t\varphi(t,x)&+\int_{\Omega}u_0(x)\varphi(0,x)dx
\\[5pt]
&=\iint_{\Omega_T} (\mu+u(t,x)) A(x)\nabla\mathcal{K}_su(t,x)\cdot\nabla \varphi(t,x) dx dt,
\end{aligned}
\label{eq:thequiaproxso13}
\end{equation}
for all test functions $\varphi\in C_c^{\infty}([0,T);C^{\infty}_{\Gamma_0}(\bar{\Omega}))$.
\end{proposition}

The proof's idea of \eqref{eq:thequiaproxso13} is to pass to the limit in \eqref{eq:thequiaproxsol} 
as $\delta \to 0^+$. First, we consider the following lemmas.
\begin{lemma}\label{Lem:weak*} Under the hypothesis of Theorem \ref{Theaprox},    there exist a subsequence of $\{u_{\delta}\}_{\delta>0}$
 such that
$$
u_\delta\rightarrow u\quad\mbox{  weakly-$\star$ in }L^{\infty}(\Omega_T),
$$ 
where $u\in L^{\infty}(\Omega_T)$.
\end{lemma}
\begin{proof}
From \eqref{eq:+ebounded}, it follows that 
$\{u_\delta\}_{\delta>0}$ is (uniformly) bounded in $L^{\infty}(\Omega_T)$.
This proves the Lemma. 
\end{proof}
\begin{lemma}\label{Lem:weak} Under the hypothesis of Theorem \ref{Theaprox},  there exist a subsequence of $\{\nabla\mathcal{K}_su_\delta\}_{\delta>0}$ and $\{u_\delta\}_{\delta>0}$ such that
$$
\begin{aligned}
\nabla\mathcal{K}_su_\delta&\rightarrow\nabla\mathcal{K}_su,\quad\mbox{  weakly in }\mathbf{L}^2(\Omega_T),
\\[5pt]
u_\delta&\rightarrow u,\quad\quad\;\; \mbox{ weakly in } L^2\big((0,T); D\big(\mathcal{L}_{\mathcal{B}}^{(1-s)/2}\big)\big),
\end{aligned}
$$
where $u\in L^2\big((0,T); D\big(\mathcal{L}_{\mathcal{B}}^{(1-s)/2}\big)\big)$.
\end{lemma}
\begin{proof}
From Proposition \ref{eq:2energystimate}, we have 
$$
    \iint_{\Omega_T} |\nabla \mathcal{K}_s u_{\delta}|^2 \,dxdt
    \leq \dfrac{C}{\mu},
$$
where $C$ is a positive constant which does not depend on $\delta$.
Therefore, the right-hand side is (uniformly) bounded in $\mathbf{L}^2(\Omega_T)$ 
w.r.t. $\delta$. Thus we obtain (along suitable subsequence) that,
$\nabla\mathcal{K}_su_\delta$ converges weakly to $\textbf{v}$ in $\mathbf{L}^2(\Omega_T)$.

\medskip
The next step is to show that $\textbf{v}= \nabla\mathcal{K}_su$ in $\mathbf{L}^2(\Omega_T)$. 
First we prove the regularity of $u$. From the equivalent norm \eqref{normequi}  we deduce that 
\begin{equation}
\begin{aligned}
   \iint_{\Omega_T} \left|\mathcal{L}_{\mathcal{B}}^{(1-s)/2} u_\delta(t,x)\right|^2dxdt &
 \leq\Lambda_2\iint_{\Omega_T} \vert\nabla \mathcal{H}_s u_{\delta}(t,x) \vert^2dxdt.
\end{aligned}\nonumber
\end{equation}
 
 On the other hand, from Corollary \ref{cor:gradlimi}, we obtain that $\nabla\mathcal{H}_s u_{\delta}$ is (uniformly) bounded in $\mathbf{L}^2(\Omega_T)$ w.r.t. $\delta$.
Thus $\{u_\delta\}$ is (uniformly) bounded in  $L^2\big((0,T); D\big(\mathcal{L}_{\mathcal{B}}^{(1-s)/2}\big)\big)$. Consequently, it is possible to select a 
subsequence, still denoted by $\{u_\delta\}$, converging weakly to $u$ in $L^2\big((0,T); D\big(\mathcal{L}_{\mathcal{B}}^{(1-s)/2}\big)\big)$, 
where we have used the uniqueness of the limit. Therefore, using again \eqref{normequi} and  the Poincare's type inequality  (Corollary \ref{poincare}), it follows that
$$
\iint_{\Omega_T} |\nabla\mathcal{K}_su(t,x)|^2 dxdt\leq\Lambda_1^{-1}\lambda_1^{-s}\iint_{\Omega_T}|\mathcal{L}_{\mathcal{B}}^{(1-s)/2}u(t,x)|^2dxdt,
$$
where $\lambda_1$ is the first eigenvalue of $\mathcal{L}$. Thus, we obtain that $\nabla\mathcal{K}_su\in \mathbf{L}^2(\Omega_T)$, and 
hence $\nabla\mathcal{K}_su_\delta$ converges weakly to 
$\nabla\mathcal{K}_su$  in $\mathbf{L}^2(\Omega_T)$.
\end{proof}
\begin{lemma}\label{Lem:strong}
Under the hypothesis of Theorem \ref{Theaprox},  there exist a subsequence of $\{u_\delta\}_{\delta>0}$ such that, 
$$
u_\delta\rightarrow u\quad\mbox{  strongly in }L^2(\Omega_T),
$$
where $u\in L^2\big((0,T); D\big(\mathcal{L}_{\mathcal{B}}^{(1-s)/2}\big)\big)$.
\end{lemma}
\begin{proof} 
Here we apply the Aubin-Lions compactness Theorem. First, from Lemma \ref{Lem:weak} we have 
$$
u_\delta\rightarrow u, \mbox{ weakly in } L^2\big((0,T); D\big(\mathcal{L}_{\mathcal{B}}^{(1-s)/2}\big)\big).
$$

On the other hand, from  Proposition \ref{eq:Fu}, \ref{eq:2energystimate} and \ref{eq:L2Hs-1}, together with the (uniform) boundedness of $\nabla\mathcal{K}_su_\delta$ in  $\mathbf{L}^2(\Omega_T)$, we have
\begin{equation}
\int^T_0 \left\Vert \partial_t u_\delta\right\Vert^2_{H^{-1}(\Omega)} dt \leq C \ ( \Vert u_0 \Vert_{\infty} + \mu).
\label{eq:L2H-14.2}
\end{equation}
One observes that, at this point $\mu > 0$ is fixed. Thus, the right-hand side of \eqref{eq:L2H-14.2} is  bounded in $L^2((0,T);H^{-1}(\Omega))$ 
w.r.t. $\delta$. Therefore, exist a subsequence, such that $\partial_t u_\delta$ converges weakly to $\partial_t u$ in 
$L^2(0,T;H^{-1}(\Omega))$. Then, applying the Aubin-Lions compactness Theorem (see \cite{Malek}, Lemma 2.48)
it follows that, $u_\delta$ converges to $u$ (along suitable subsequence) strongly in 
$L^2(\Omega_T)$ as $\delta$ goes to zero. 
\end{proof}

\begin{proof}[\textbf{Proof of Proposition \ref{Thprinc}}]
The idea of the proof of \eqref{eq:thequiaproxso13} is to pass to the limit in 
\eqref{eq:thequiaproxsol} as $\delta \to 0^+$.
From Lemma \ref{Lem:weak*} is enough to pass to the limit in the first integral in the left hand side of \eqref{eq:thequiaproxsol}. We can proceed in a similar way as before for the sequence $u_{0,\delta}$. 

\medskip
On the other hand, by Corollary \ref{cor:gradlimi} and H\"older inequality, we have that the second integral in the left hand side of \eqref{eq:thequiaproxsol} is zero, given that $A\in L^{\infty}(\Omega)$ and $\phi\in L^2(\Omega)$.

\medskip
Now, we study the convergence of the integral in right hand side of \eqref{eq:thequiaproxsol}. First, since $A(x)$ is symmetric,  it is sufficient to show $(\mu+u_\delta)\nabla\mathcal{K}_su_\delta$ converges weakly in $\mathbf{L}^2(\Omega_T)$. Indeed, by Lemma \ref{Lem:weak} and \ref{Lem:strong}, we obtain that
$(\mu+u_\delta) \nabla\mathcal{K}_su_\delta$ converges weakly to $(\mu+u)\nabla\mathcal{K}_su$ as $\delta\to 0^+$.
Hence, the equality \eqref{eq:thequiaproxso13} follows. 
\end{proof}
\begin{corollary}\label{corprinc}
Let $u$ be the function given by Proposition \ref{Thprinc}, then it satisfies:
\item[$(1)$] For almost all  $(t,x) \in \Omega_T$
\begin{equation}
\label{eq:bounded4.2}
   0\leq u(t)+\mu \leq  \Vert u_0 \Vert_{\infty}, \quad \text{and} 
\end{equation}
\begin{equation}
\label{eq:consevationofmass4.2}
    \int_{\Omega}u(x,t)dx = \int_{\Omega}u_0(x)dx.
\end{equation}

\item[$(2)$] First energy estimate: For 
$\eta(\lambda):= (\lambda+\mu) \log (1+(\lambda/\mu))-\lambda$, $(\lambda\geq 0)$, and almost all $t\in (0,T)$, 
\begin{equation}
\int_{\Omega}\eta(u(t)) dx+\Lambda_1\int_0^t\int_{\Omega}\vert \nabla \mathcal{H}_su\vert^2 \ dxdt'
\leq\int_{\Omega}\eta(u_0) \ dx.
\label{energyinequa4.2}
\end{equation} 

\item[$(3)$] Second energy estimate:  For almost all $0< t_1< t_2< T$,
\begin{equation}
\frac{1}{2}\int_{\Omega}\vert \mathcal{H}_s u(t_2) \vert^2 dx + \Lambda_1\int_{t_1}^{t_2}\int_{\Omega} (\mu+u) 
\vert \nabla \mathcal{K}_su \vert^2 \,dx\,dt\leq\frac{1}{2}\int_{\Omega}\vert \mathcal{H}_s u(t_1) \vert^2 dx.\label{eq:2energystimate4.2}
\end{equation}

\item[$(4)$] For each $v \in H^1_{\Gamma_0}(\Omega)$,
\begin{equation}
\int^T_0 \langle \partial_t u, v\rangle dt= \iint_{\Omega_T} (\mu+u)A(x)\nabla \mathcal{K}_su\cdot \nabla v \,dx dt, 
\label{estimateH-14.2}
\end{equation}
where $\langle\cdot,\cdot\rangle$ denotes the pairing between $\left(H^1_{\Gamma_0}(\Omega)\right)^*$ and $H^1_{\Gamma_0}(\Omega)$.
\end{corollary}
\begin{proof}
1. To show \eqref{eq:bounded4.2}. Recall that $u_\delta$ converges strong to $u$ in $L^2(\Omega_T)$ and therefore (for a subsequence) 
$u_\delta$ converges a.e. to $u$ in $\Omega_T$, then passing the limit in  \eqref{eq:+ebounded} as $\delta \to 0^+$, we obtain the \eqref{eq:bounded4.2}. Assertion \eqref{eq:consevationofmass4.2} is obtained by 
\eqref{eq:consermassaprox} together with 
\medskip
 Dominated Convergence Theorem.

\medskip
2. To prove the first energy estimate \eqref{energyinequa4.2}, we pass to the limit in \eqref{eq:firstenergy} as $\delta\to 0^+$. 
Due to $u_\delta$ converges almost everywhere to $u$ in $\Omega_T$, and $\eta$ is a continuous function,
it follows that $\eta(u_\delta)$ converges 
almost everywhere to $\eta(u)$ in $\Omega_T$. Moreover, $u_\delta$ is bounded in $L^{\infty}(\Omega_T)$ w.r.t. $\delta$, 
then for almost all $t\in (0,T)$
$$
    \lim_{\delta\to 0^+}\int_{\Omega}\eta(u_\delta(t)) \ dx= \int_{\Omega}\eta(u(t)) \ dx,
$$ 
where we have used  the Dominated Convergence Theorem. 
We can proceed in a similar way as before for the sequence $u_{0,\delta}$.

\medskip
On the other hand, using the  idea of the proof of Lemma \ref{Lem:weak} 
it is possible to show that (for a subsequence) $\nabla\mathcal{H}_su_\delta$ converges weakly to $\nabla\mathcal{H}_su$ in 
$\mathbf{L}^2(\Omega_T)$. Then, we have
$$
\int_0^t\int_{\Omega} |\nabla \mathcal{H}_su|^2 \,dxdt' \leq \liminf_{\delta\to 0^+} \int_0^t\int_{\Omega} |\nabla \mathcal{H}_su_\delta|^2 \,dxdt'
$$
for almost all $t \in (0,T)$.  Also observe that the second integral in the left hand side of \eqref{eq:firstenergy} is positive, hence we throw it out. 
Therefore passing to the limit in \eqref{eq:firstenergy} as $\delta$ tends to zero, we obtain the assertion.

\medskip
3. To show the second energy estimate \eqref{eq:2energystimate4.2}, we pass to the limit in 
\eqref{eq:secondenergy} as $\delta$ goes to zero. First, we have to study the convergence of 
 each integral in \eqref{eq:2energystimate}. One notes that, due to the continuity in $L^2(\Omega_T)$ and Lemma \ref{Lem:strong}, 
it follows that $\mathcal{H}_su_\delta$ strongly converges to $\mathcal{H}_su$ in $L^2(\Omega_T)$. 
Consequently, it is possible to select a subsequence, still denoted by $\mathcal{H}_su_\delta(t)$  such that,
for almost all $t\in(0,T)$
$$
\lim_{\delta\to 0^+}\int_{\Omega}|\mathcal{H}_su_\delta(t,x)|^2\,dx=\int_{\Omega}|\mathcal{H}_su(t,x)|^2\,dx.
$$
On the other hand, since
second integral in the left hand side of \eqref{eq:secondenergy} is positive for all $\delta>0$, hence we throw it out. Finally, the convergence of the 
third integral follows from Lemma \ref{Lem:weak} and \ref{Lem:strong}. Then, passing 
to the limit in \eqref{eq:secondenergy} as $\delta \to 0^+$, 
we obtain \eqref{eq:2energystimate4.2}.

\medskip
4. Assertion \eqref{estimateH-14.2} follows by similar ideas, so we pass to the limit in \eqref{eq:3estimativa} as $\delta\to 0^+$, and
the proof is concluded.

\end{proof}
\begin{remark}
\label{remark4.2}
The function $u$ obtained above depends on the fixed parameter $\mu$.  
For each $\mu> 0$, we write from now on $u_\mu$ instead of $u$. 
\end{remark}

\subsubsection{Limit transition $\mu \to 0^+$}
Here, we prove the existence of weak solutions for the initial mixed-boundary value problem \eqref{FTPME}. To show that  
we consider the sequence $\{u_\mu\}_{\mu>0}$, obtained in  Proposition \ref{Thprinc}, 
which satisfies  Corollary \ref{corprinc} for each $\mu> 0$, \eqref{eq:thequiaproxso13}--\eqref{estimateH-14.2}.

\begin{proof}[\textbf{Proof of Theorem \ref{Thprincipal}}]

To show the existence of solution we  pass to the limit in \eqref{eq:thequiaproxso13} as $\mu\to 0^+$. From \eqref{eq:bounded4.2} and $\mu\in(0,1)$, we see that 
$\{u_\mu\}_{\mu>0}$ is (uniformly) bounded in $L^{\infty}(\Omega_T)$ w.r.t $\mu$.
Hence, it is possible to select a subsequence, still denoted by $\{u_\mu\}$,
converging weakly-$\star$ to $u$ in $L^{\infty}(\Omega_T)$, 
which is enough to pass to the limit in the first integral in the left hand side of \eqref{eq:thequiaproxso13}. 

\medskip
Now, we study the convergence of the integral in right hand side of \eqref{eq:thequiaproxso13}. First, since $A(x)$ is symmetric,  
it is sufficient to show $(\mu+u_\mu)\nabla\mathcal{K}_su_\mu$ converges weakly in $\mathbf{L}^2(\Omega_T) $. 
On the other hand, we recall that, for each $\lambda \geq 0$, 
$$
\begin{aligned}
\eta(\lambda)&= (\lambda+\mu)\log (1+\lambda/\mu)- \lambda, 
\\[5pt]
&=(\lambda+\mu)\log (\lambda+\mu)-(\lambda+\mu)\log \mu-\lambda.
\end{aligned}
$$
Then, from \eqref{eq:consevationofmass4.2} and \eqref{energyinequa4.2} we obtain for almost all $t\in (0,T)$
\begin{equation}
\begin{aligned}
    \Lambda_1\int_0^t\int_\Omega |\nabla \mathcal{H}_su_\mu|^2 \ dxdt
    &+\int_\Omega (u_\mu(t)+\mu)\log(u_\mu(t)+\mu) \ dx
\\[5pt]
&\leq \int_{\Omega} (u_0 +\mu) \log(u_0+\mu) \ dx.
\end{aligned}\label{eq:deslog}
\end{equation}

\medskip
Since $f= f^+-f^-$, where $f^\pm= \max\{\pm f, 0\}$, it follows from \eqref{eq:deslog} that
$$
\begin{aligned}
   \Lambda_1\int_0^t\int_\Omega &|\nabla \mathcal{H}_su_\mu|^2 \ dxdt
   +\int_\Omega (u_\mu(t)+\mu) \log^{+}(u_\mu(t)+\mu) \ dx
\\[5pt]
\leq & \int_{\Omega} (u_0+\mu) \log(u_0+\mu) \ dx+\int_\Omega (u_\mu(t)+\mu) \log^-(u_\mu(t)+\mu)dx.
\end{aligned}
$$

We observe that, the right hand side of the above inequality is bounded w.r.t. 
$\mu$ (small enough), because $u_\mu$ is bounded in $L^{\infty}(\Omega_T)$ w.r.t. $\mu$, 
and
$$
\int_\Omega (u_\mu(t)+\mu)\log^-(u_\mu(t)+\mu)dx,
$$
is bounded w.r.t. $\mu$  (small enough).
Consequently, we have that
$\nabla \mathcal{H}_su_{\mu}$ is (uniformly) bounded in $L^2(\Omega_T)$.

\medskip
On the other hand, using \eqref{normequi} and the Poincar\'e inequality (Corollary \ref{poincare}), we obtain that 
\begin{equation}
\begin{aligned}
   \iint_{\Omega_T} \left|\nabla \mathcal{K}_su_\mu(t,x)\right|^2dxdt&\leq \Lambda_1^{-1} \iint_{\Omega_T} \left|\mathcal{L}_{\mathcal{B}}^{1/2-s}u_\mu(t,x)\right|^2dxdt
   \\[5pt]
& \leq\Lambda_1^{-1}\lambda_1^{-s}\iint_{\Omega_T} \vert\mathcal{L}_{\mathcal{B}}^{1/2-s/2}u_\mu(t,x)\vert^2dxdt
\\[5pt]
&\leq \Lambda_1^{-1}\lambda_1^{-s}\Lambda_2\iint_{\Omega_T} \vert\nabla\mathcal{H}_su_\mu(t,x)\vert^2dxdt.
\end{aligned}\nonumber
\end{equation}
Therefore, $\nabla \mathcal{K}_su_\mu$ is (uniformly) bounded in $L^2(\Omega_T)$ w.r.t. $\mu> 0$, and thus we obtain (along suitable subsequence) 
that $\nabla\mathcal{K}_su_\mu$ converges weakly to $\textbf{v}$ in $\mathbf{L}^2(\Omega_T)$.
It remains to show that $\textbf{v}=\nabla\mathcal{K}_su$. 
Moreover, applying the same ideas as in the proof of the Proposition \ref{Thprinc}, it is possible to select a subsequence, 
still denoted by $\{u_\mu\}$, converging weakly to $u$ in $L^2\left(0,T;D\big( \mathcal{L}_{\mathcal{B}}^{(1-s)/2}\big)\right)$, 
such that 
$$
   \text{$\textbf{v}=\nabla\mathcal{K}_su$ in $\mathbf{L}^2(\Omega_T)$.} 
$$ 
Hence $\nabla\mathcal{K}_su_\delta$ 
converges weakly to $\nabla\mathcal{K}_su$  in $\mathbf{L}^2(\Omega_T)$.

\medskip
Now, we prove strong convergence for $\{u_\mu\}_{\mu>0}$ in $L^2(\Omega_T)$. To show that, we apply again the Aubin-Lions 
compactness Theorem. 
Since the coefficient of the matrix $A(x)$ are in $ C(\overline{\Omega}) \cap C^{0,1}_{\loc}(\Omega)$, 
together with the boundedness of $\nabla\mathcal{K}_su_\mu$ in $\mathbf{L}^2(\Omega_T)$, 
and the uniform limitation of $u_\mu$, we have from \eqref{estimateH-14.2} that
\begin{equation}
\int^T_0 \left\Vert \partial_t u_\mu \right\Vert^2_{H^{-1}(\Omega)} dt \leq C,
\label{eq:L2H-14.3}
\end{equation}
where $C$ is a positive constant which does not depend on $\mu$. Then,
passing to a subsequence (still denoted by $\{u_\mu\}$), we obtain that 
$$
  \text{$\partial_t u_\mu$ converges weakly to
 $\partial_t u$ in $L^2(0,T;H^{-1}(\Omega))$.}
$$
Applying the Aubin-Lions compactness Theorem, it follows that $u_\mu$ converges strongly to $u$ (along suitable sequence) in $L^2(\Omega_T)$.
Consequently, we obtain that $(\mu+u_\mu)\nabla\mathcal{K}_su_\mu$ converges weakly to $u \ \nabla\mathcal{K}_su$ as $\mu\to 0^+$. Then, we are ready to pass to the limit in \eqref{eq:thequiaproxso13} as $\mu\to 0^+$ to get 
\begin{equation}
    \iint_{\Omega_T} u(t,x) \big( \partial_t\varphi(t,x) -  A(x)\nabla \mathcal{K}_s(u(t,x)) \cdot \nabla \varphi(t,x) \big)  \ dxdt
   +\int_{\Omega}u_0(x)\varphi(0,x)dx=0,\nonumber
\end{equation}
for all $\vp \in C^{\infty}_c([0,T);C^{\infty}_{\Gamma_0}(\bar{\Omega}) )$. 
\end{proof}
\begin{corollary}
The  solution $u$ of the initial mixed-boundary value problem \eqref{FTPME} given by Theorem \ref{Thprincipal},
satisfies: 
\begin{enumerate}
\item[$(1)$]
For almost  all $t \in (0,T)$, we have
\begin{equation}
\Vert u(t)\Vert_{\infty} \leq  \Vert u_0 \Vert_{\infty}, \quad \text{and} \label{eq:bounded4.3}
\end{equation}
\begin{equation}
\int_{\Omega}u(x,t)dx = \int_{\Omega}u_0(x) \ dx. \label{eq:consevationofmass4.3}
\end{equation}

\item[$(2)$] First energy estimate: For almost all  $t \in (0,T)$, 
\begin{equation}
    \Lambda_1\int^t_0\int_{\Omega}\vert \nabla \mathcal{H}_su \vert^2 \ dxdt'
   +\int_{\Omega} u(t) \log(u(t)) \ dx \leq \int_{\Omega} u_0 \log(u_0) \ dx.\label{energyinequa4.3}
\end{equation}
\item[$(3)$] Second energy estimate: For almost all $0<t_1<t_2<T$,
\begin{equation}
\frac{1}{2}\int_{\Omega}\vert \mathcal{H}_s u(t_2) \vert^2 \ dx  + 
\Lambda_1\int_{t_1}^{t_2} \int_{\Omega} u\vert \nabla \mathcal{K}_su\vert^2 \,dx\,dt \leq\frac{1}{2}\int_{\Omega}\vert \mathcal{H}_s u(t_1) \vert^2 dx.
\label{eq:2energystimate4.3}
\end{equation}
\end{enumerate}
\end{corollary}

\begin{proof} In order
to show \eqref{eq:bounded4.3}-\eqref{eq:2energystimate4.3},
we may follow similar lines as in the proof of Corollary \ref{corprinc}. Therefore, we omit them here.
\end{proof}

\section{Appendix}

Let us fix here some notation and background used in this paper, we first consider the notion 
of $C^1$-(admissible) deformations, which is used to give the correct notion of traces. 
One can refer to \cite{NPS}.

\begin{definition}
\label{AdmDefo}
Let $\Omega \subset \R^{n}$ be an open set.
A $C^1$-map $\Psi:[0,1] \times \Gamma\to \ol{\Omega}$ is said a $C^1$ 
admissible deformation,
when it satisfies the following conditions:
\begin{enumerate}
\item[$(1)$] For all $r\in\Gamma$, $\Psi(0,r)=r$.

\item [$(2)$] The derivative of the map $[0,1] \ni \tau \mapsto \Psi(\tau, r)$
at $\tau= 0$ is not orthogonal
to $\nu(r)$, for each $r\in \Gamma$.
\end{enumerate}
\end{definition}
Moreover, for each $\tau\in[0,1]$, we denote:
$\Psi_\tau$ the mapping from $\Gamma$ to $\ol{\Omega}$, given by $\Psi_\tau(r):=\Psi(\tau,r)$;
$\nu_\tau$ the unit outward normal field in $\Psi_\tau(\Gamma)$. In particular,
$\nu_0(x)=\nu(x)$ is the unit outward normal field in $\Gamma$.

\medskip
It must be recognized that domains with $C^2$ boundaries 
always have $C^1$ admissible deformations. Indeed, 
it is enough to take $\Psi(\tau,r)= r - \epsilon \tau \nu(r)$ for
sufficiently small $\epsilon > 0$. 

\medskip
Now, we define a level set function $h$ associated with the deformation $\Psi_\tau$.
For $\delta> 0$ sufficiently small we define
$$
   h(x):= \left \{ 
         \begin{aligned}
           \min\{\tau,\delta\}, &\quad \text{if $x \in \Omega$},
            \\[5pt]
           - \min\{\tau,\delta\},& \quad \text{if $x \in \R^n \setminus \Omega$},           
         \end{aligned}\right.
$$
which is Lipschitz continuous in $\R^n$, and $C^1$ on the closure of $\left\lbrace x \in \R^n: |h(x)| < \delta \right\rbrace$, moreover
\begin{equation}
|\nabla h(x)|=\left\{\begin{array}{ll}
1 & \text { for } 0 \leq h(x)<\delta ,
\\[5pt]
0 & \text { for } h(x)=\delta.
\end{array}\right.\label{eq:nablah}
\end{equation}

\begin{lemma}
\label{Lemma:xik}
Let $\Omega \subset \R^n$ be an open bounded domain with $C^2$ boundary. For each $k \in \Bbb{N}$, and all $x\in \R^n$, consider
\begin{equation}
\xi_{k}(x):= 1-{\exp} \left( -k \ h(x) \right).\label{eq:xi}
\end{equation}
Then,  the sequence $\{\xi_{k}\}$ satisfies
\begin{equation}
 \lim_{k\to +\infty}\int_{\Omega}|1-\xi_k|^2dx=0,
  \quad \text{and} \quad
  \lim_{k\to +\infty}\int_{\Omega}|\nabla\xi_k|^2dx=0.
  \label{eq:limvarep}
\end{equation}
\end{lemma}
\begin{proof}
For more details see M\'alek, Necas, Rokyta and Ruzicka \cite{Malek}, p. 129.
\end{proof}

Last but not least, let us consider the following approximating sequences. 
Choose a non-negative function $\gamma \in C^1_c(\R)$, with support contained in $[0,1]$, 
such that, $\int \gamma(t) dt= 1$. Then, we consider the sequences $\{\delta_j\}_{j\in \N}$,
and $\{H_j\}_{j\in \N}$, defined by
$$
   \delta_j(t):= j \ \gamma(j t), \quad H_j(t):= \int_0^t \delta_j(s) \ ds.
$$
Thus, $H'_j(t)= \delta_j(t)$, and clearly the sequence $\delta_j(t)$  
converges as $j \to \infty$ to the Dirac $\delta$-measure in $\mathcal{D}'(\R)$,
while the sequence $H_j(t)$  converges pointwise to the Heaviside function 
$$
   H(t)= \left \{
   \begin{aligned}
   1, &\quad \text{if $t \geq 0$},
   \\
   0, &\quad \text{if $t< 0$}.
   \end{aligned}
   \right.  
$$

\medskip
To finish this section, we show the existence and uniqueness of $u_{\mu,\delta}$
for the approximate parabolic problem \eqref{eq:aproxequation1}--\eqref{eq:aproxequation4}. 
To this end, 
we first apply the Banach Fixed Point Theorem to prove the local in time existence of solution,
and thus applying a contradiction argument we extend it to be global in time. Since \eqref{eq:aproxequation1} is a fractional non-standard 
parabolic equation, we present the important details and omit the usual ones. 

\begin{theorem}
\label{EXUNIQAPP}
Let $u_{0\delta}$ be a non-negative 
regularized initial data. Then the problem \eqref{eq:aproxequation1}--\eqref{eq:aproxequation4} admits a unique strong solution 
$$
\begin{aligned}
u_{\mu,\delta}&\in C([0,T); H^1_{\Gamma_0}(\Omega))\cap L^2((0,T); H^2(\Omega^{\prime}))  \cap L^\infty(\Omega_T), 
\\[5pt]
\partial_t u_{\mu,\delta} &\in L^2\left(\Omega_{T}\right),
\end{aligned}
$$
for each $\Omega^{\prime}$ compactly contained in $\Omega$.
\end{theorem}
 \begin{proof}
The proof will be divided into four steps. 

\medskip 
1. First, for each  $\tilde{u}\in L^{\infty}(\Omega_T) \cap L^2(0,T;D(\mathcal{L}_{\mathcal{B}}^{1-s}))$, the following problem
\begin{equation}
\label{Frac_para}
\begin{cases}
     \partial_t u_{\mu,\delta}
     -\delta \ \dive(A(x)\nabla u_{\mu,\delta})=
     \dive(\mathbf{q}_{\mu}(x,\tilde{u})) 
     &\mbox{in}\;\Omega_T,
     \\[5pt]
{u_{\mu,\delta}}=u_{0\delta} &\mbox{in} \{t=0\} \times \Omega,
\\[5pt]
 u_{\mu,\delta}= 0 &\mbox{on}\; (0,T) \times \Gamma_0,
 \\[5pt]
\delta A\nabla u_{\mu,\delta}\cdot \nu =-\mathbf{q}_{\mu}(x,\tilde{u})\cdot \nu&\mbox{on}\; (0,T) \times \Gamma_1,
\end{cases}
\end{equation}
has a unique weak solution 
$$
 u_{\mu,\delta} \in L^2\big((0,T);H^1_{\Gamma_0}(\Omega)\big)\cap C\big([0,T);L^2(\Omega)\big) \cap L^\infty(\Omega_T). 
$$    
Indeed, since $\tilde{u}\in L^{\infty}(\Omega_T) \cap L^2(0,T;D(\mathcal{L}_{\mathcal{B}}^{1-s}))$, 
it follows that 
$$
    \mathbf{q}_{\mu}(x,\tilde{u}) \in L^{2}((0,T); H^1_{\Gamma_0}(\Omega)).
$$ 
Then applying the parabolic theory, see Theorem 11.8 in Chipot \cite{MCHIPOT}, (also Chipot, Rougirel \cite{MCHIPOTROUGIREL}),
there exists a unique weak solution 
$$
u_{\mu,\delta} \in L^2\big((0,T); H_{\Gamma_0}^1(\Omega)\big)\cap C\big([0,T); L^2(\Omega)\big) \cap L^\infty(\Omega_T) 
$$ 
of the problem \eqref{Frac_para}. 

\medskip
2. Now, we show the local in time existence of solution 
to \eqref{eq:aproxequation1}--\eqref{eq:aproxequation4}. 
To prove that, we define the following map 
$$
\begin{aligned}
u_{\mu,\delta}(t,x)= \mathcal{T}(\tilde{u})(t,x)&:= \int_{\Omega}K(t,x,y) \ u_{0,\delta}(y)dy
\\[5pt]
&+ \int_0^t \!\! \int_{\Omega}(\tilde{u}(t',y)+\mu)  \nabla_y K(t-t',x,y) \cdot \nabla\mathcal{K}_s\tilde{u}(t',y) dy dt,
\end{aligned}
$$
where $K(t,x,y)$, $(x,y \in \Omega)$, is the heat kernel of the operator $\mathcal{L}u \!=-\dive(A(\cdot)\nabla u)$ with 
mixed Dirichlet-Neumann boundary data, see \cite{DAVIES}. Moreover,
for $t> 0$ sufficiently small, it is not difficult to show that $\mathcal{T}$ is a contraction. 
Then, applying the Banach Fixed Point Theorem, there exists a unique local in time weak solution
$$u_{\mu,\delta} \in L^2\big((0,T_M);H_{\Gamma_0}^1(\Omega)\big)\cap C\big([0,T_M);L^2(\Omega)\big) \cap L^\infty(\Omega_{T_M}),$$
where $T_M$ denotes the maximal time of existence. 

\medskip
3. We claim that the local solution 
$u_{\mu,\delta}$ satisfies  
\begin{equation}
\begin{aligned}
u_{\mu,\delta}&\in C\left([0,T_M);H^1_{\Gamma_0}(\Omega)\right)\cap L^2\left((0,T_M);H^2(\Omega^{\prime})\right) \cap L^\infty(\Omega_{T_M}), 
\\[5pt]
\partial_t u_{\mu,\delta} &\in L^2\left(\Omega_{T_M}\right).
\end{aligned}
\label{eq:regul_ape}
\end{equation}

Indeed, since $u_{\mu,\delta} \in L^2\big((0,T_M);H_{\Gamma_0}^1(\Omega)\big)\cap C\big([0,T_M);L^2(\Omega)\big) \cap L^\infty(\Omega_{T_M})$, we have 
$$
\dive\left((u_{\mu,\delta}+\mu) A(x) \nabla\mathcal{K}_{s}u_{\mu,\delta} \right)\,\in\, L^2\left((0,T_M);L^2(\Omega)\right).
$$ 
Therefore, from equation \eqref{eq:aproxequation1} and the standard parabolic regularity theory (see \cite{ACM}), we obtain \eqref{eq:regul_ape}.
Consequently, $u_{\mu,\delta}$ satisfies the partial differential equation \eqref{eq:aproxequation1}
in the strong sense, that is, for almost all $(t,x) \in (0, T_M) \times \Omega^{\prime}$. 

%

\medskip
4. Finally, we claim that $T_M= T$, for any $T> 0$. Conversely, let us 
suppose that, $T_M< T$. Then, there exists an increasing sequence $\{t_j\}_{j=1}^\infty$, such that,
$t_j\to T_M^{-}$ as $j\to\infty$ and
\begin{equation}
\label{BLOWUP}
\lim_{j\to\infty}\|u_{\mu,\delta}(t_j,\cdot)\|_{L^{\infty}(\Omega)}=+\infty.
\end{equation}
Although, due to a similar proof given to \eqref{eq:+ebounded},
we may show that
$$
   0\leq  u_{\mu,\delta}(t,x)+\mu\leq\Vert u_{0\delta}\Vert_{L^{\infty}(\Omega)},
$$
for each $t \in (0,T_M)$ and almost all $x \in \Omega$,
which contradicts \eqref{BLOWUP}. 
\end{proof}
 \section*{Acknowledgements}

We are grateful to the anonymous referee for carefully reading the previous 
version of this manuscript and for valuable comments, which helped us to 
improve the exposition of this paper. 

\medskip
Conflict of Interest: Author Wladimir Neves has received research grants from CNPq
through the grant  308064/2019-4, and also by FAPERJ 
(Cientista do Nosso Estado) through the grant E-26/201.139/2021.



\begin{thebibliography}{99}

\bibitem{HA} {\sc H. Amann}, {\em Nonhomogeneous linear and quasilinear elliptic and parabolic boundary value problems, Function Spaces, Differential
Operators and Nonlinear Analysis} (Eds. Schmeisser, Triebel), Teubner Texte zur Mathematik 133 (1993), pp. 9--126.

\bibitem{ACM} {\sc L. Ambrosio, A. Carlotto, A. Massaccesi},
{\em Lectures on Elliptic Partial Differential Equations}.
Springer, 2018.

\bibitem{BIK} {\sc P. Biler, C. Imbert, G. Karch }, {\em The nonlocal porous medium equation: Barenblatt profiles and other weak solution}, Arch. Ration. Mech. Anal. 215(2015), 497529. 

\bibitem{MSireVazquez} {\sc M. Bonforte, Y. Sire, J. L. V\'azquez}, {\em Existence, Uniqueness and Asymptotic behaviour
for fractional porous medium equations on bounded domains}, Manuscript submitted to AIMS' Journals.

\bibitem{CSV} {\sc L. Caffarelli, F. Soria, J. L. Vazquez}, {\em Regularity of solution of the fractional porous medium flow} J. Eur. Math. Soc. (JEMS) 155 (2013), 17011746.

\bibitem{Caffa} {\sc L. Caffarelli, J. L. Vazquez}, {\em Nonlinear Porous Medium Flow with Fractional Potential Pressure}, Arch. Rational Mech. Anal. 202 (2011), 537--565.

\bibitem{CV2} {\sc L. Caffarelli, J. L. Vazquez}, {\em Regularity of solution of the fractional porous medium flow in the exponent $1/2$}, St Petersburg Math Journal 27(2016), no 3, 437460.

\bibitem{LCPRS} {\sc L. Caffarelli, P. R. Stinga}, {\em Fractional elliptic equations, Caccioppoli estimates and regularity},
Annales de l'Institut Henri Poincar\'e C, Analyse Non Lin\'eaire, 33 (2016), 767--807.

\bibitem{MCHIPOT}
{\sc M. Chipot}, {\em Elements of Nonlinear Analysis}, 
Birkh\"auser, 2000.

\bibitem{MCHIPOTROUGIREL}
{\sc M. Chipot, A. Rougirel}, {\em On some class of problems with nonlocal source and boundary flux},
Advances in Differential Equations, Vol 6,  No 9 (2001), 1025--1048.

\bibitem{DAVIES}
{\sc E. B. Davies}, {\it Heat kernels and spectral theory}.
Cambridge University Press, 1989. 

\bibitem{Grubb} {\sc G. Grubb}, {\em Regularity of spectral fractional Dirichlet and Newmann problems}, Math. Nachr. 289 (2016) 7, 831--844.

\bibitem{Hua} {\sc G. Huaroto, W. Neves}, {\em Initial-boundary value problem for a fractional type degenerate heat equation}, Mathematical Models and Methods
in Applied Sciences, Vol. 28, No. 6 (2018) 1199-1231.

\bibitem{GHWN3} {\sc G. Huaroto, W. Neves}, {\em Solvability of the Fractional Hyperbolic Keller-Segel System},
Submitted. 

\bibitem{I} {\sc C. Imbert}, {\em Finite speed of propagation for a non-local porous medium equation}, Colloq. Math. 143(2) (2016), 149157.


\bibitem{LionsMagenes} {\sc J. L. Lions, E. Magenes}, {\em Problemes aux limites non-homogenes et application}, V.1, Dunod, Paris, 1968.

\bibitem{LMS} {\sc S. Lisini, E. Mainini, A. Segatti.}, {\em A gradient flow approach to the porous medium equation with fractional pressure}, Arch. Rat. Mech. Anal. 227, 567-606. 

\bibitem{Malek} {\sc J. M\'alek, J. Necas, M. Rokyta, M. Ruzicka}, {\em Weak and Measure-valued solutions to evolutionary PDEs}, Chapman and  Hall, London, 1996.

\bibitem{NPS} {\sc W. Neves, E. Panov, J. Silva}, {\em Strong Traces for Conservation Laws with General Nonautonomous Flux}
SIAM Journal on Mathematical Analysis 50 (6), 6049-6081

\bibitem{SSJV} {\sc S. Serfaty, J. L. V\'azquez}, 
{\em A mean field equation as limit of nonlinear diffusions with fractional Laplacian operators},
Calculus of Variations and Partial Differential Equations volume 49 (2014), 1091--1120.

\bibitem{STV2} {\sc D. Stan, F. del Teso, J. L. Vazquez}, 
{\em Finite and infinite speed of propagation for porous medium equation with nonlocal pressure}, J. Diff. Eq., 260 1154-1199, 2016.

\bibitem{STV3} {\sc D. Stan, F. del Teso, J. L. Vazquez}, {\em Existence of weak solution solution for a general porous medium equation with nonlocal pressure}, Arch. Rat. Mech. Anal., First online 09 February 2019 in press arXiv: 1609.05139.

\bibitem{STV5} {\sc D. Stan, F. del Teso, J. L. Vazquez}, {\em Transformations of self-similar solutions for porous medium equations of fractional type}, Nolinear Anal., 119(2015), 62-73.

\end{thebibliography}
\end{document}